%% file: main.tex
\documentclass[USenglish]{article}

\usepackage[USenglish]{babel}
\usepackage{color}
\usepackage{amsmath,amsthm,amssymb,amsfonts}
\usepackage[margin=2.5 cm]{geometry}
\usepackage[unicode,colorlinks=true,linkcolor=blue,citecolor=green,urlcolor=black,pagebackref=false]{hyperref}
\usepackage{stmaryrd}
\usepackage{mathabx}
\usepackage{custom_macros_general, custom_macros_spaceops, custom_macros_funcs, custom_macros_freqops, custom_macros_deprec}

\let\oldmu\mu
\let\oldeps\varepsilon
\let\oldnu\nu
\let\oldlambda\lambda
\let\oldzeta\zeta

\renewcommand{\mu}{\boldsymbol{\oldmu}}
\renewcommand{\varepsilon}{\boldsymbol{\oldeps}}
\renewcommand{\nu}{\boldsymbol{\oldnu}}
\renewcommand{\lambda}{\boldsymbol{\oldlambda}}
\renewcommand{\zeta}{\boldsymbol{\oldzeta}}

\numberwithin{equation}{section}

\DeclareMathAlphabet\mathbfcal{OMS}{cmsy}{b}{n}
\makeatletter
\let\@fnsymbol\@arabic
\makeatother

\newcommand{\revision}[1]{{ #1}}

\title{Regularity of vector fields with piecewise regular curl and divergence}
\author{Jens Markus Melenk\thanks{(melenk@tuwien.ac.at), Institut f\"{u}r Analysis und Scientific Computing, Technische Universit\"{a}t Wien, Wiedner Hauptstrasse 8-10, A--1040 Wien, Austria.}
\and 
David W\"{o}rg\"{o}tter\thanks{(david.woergoetter@tuwien.ac.at), Institut f\"{u}r Analysis und Scientific Computing, Technische Universit\"{a}t Wien, Wiedner Hauptstrasse 8-10, A--1040 Wien, Austria.}}
\date{\today}

\begin{document}

\maketitle

\begin{abstract}
		\revision{We consider a bounded Lipschitz domain $\Omega\subseteq\R^3$ with sufficiently smooth boundary and prove piecewise Sobolev regularity of vector fields that have piecewise regular curl and divergence, but may be discontinuous across mutually disjoint and sufficiently smooth surfaces inside of $\Omega$.} The main idea behind our approach is to employ recently developed parametrices for the curl-operator and the regularity theory of Poisson transmission problems. We conclude our work by applying our findings to the heterogeneous time-harmonic Maxwell equations with either a) impedance, b) natural or c) essential boundary conditions and providing wavenumber-explicit piecewise regularity estimates for these equations.

\end{abstract}
\input{intro.tex}
	\input{notation.tex}

    \input{diffops.tex}

	\input{regularity.tex}

	\bibliographystyle{amsplain}
	\bibliography{bibliog}
\end{document}

%% file: intro.tex
\section{Introduction}

For $\ell\in\N_0$ we consider a bounded $\CM{\ell+2}$-domain \revision{$\Omega\subseteq\R^3$} with boundary $\Gamma:=\partial\Omega$ and a vector field $\solv$ which has piecewise $\VSHM{\ell}$-regular curl and $\SHM{\ell}$-regular divergence, but may jump across certain $\CM{\ell+2}$-regular and mutually disjoint surfaces $\interf_1,\ldots,\interf_r$ in the interior of $\Omega$. We show that if $\solv$ satisfies certain transmission conditions across these surfaces of discontinuity and boundary conditions on $\Gamma$, then $\solv$ is already piecewise $\VSHM{\ell+1}$-regular, and the $\VSHM{\ell+1}$-norm of $\solv$ can be controlled by its piecewise curl, its piecewise divergence and the boundary data.

\medskip

Estimates of this kind were first provided by Weber \cite{RegularityWeber} and very recently and independently of our own research, Chaumont-Frelet, Galkowski \& Spence \cite{MaxwellSpence} improved those results. In both works, the core argument is a difference quotient technique, which leads to rather technical proofs. In contrast to this, our approach is based on parametrices for the curl-operator established in \cite{CurlInverse,MaxwellImpedanceMelenk} and the regularity of Poisson transmission problems.

This new approach does not only yield $\VSHM{\ell+1}$-regularity but also slightly sharpens the results from \cite{RegularityWeber, MaxwellSpence} and provides Helmholtz-type decompositions of $\solv$, thus giving new insights on the general structure of vector fields with piecewise regular curl and divergence. In addition to that, we allow for non-homogeneous boundary conditions of $\solv$ on $\Gamma$, whereas the existing works \cite{RegularityWeber, MaxwellSpence} restrict themselves to homogeneous tangential- or homogeneous normal traces of $\solv$.  
\medskip

Our primary reason for studying vector fields with (piecewise) regular curl and divergence is their role in the solution theory of Maxwell's equations. To illustrate this, we apply our findings to the heterogeneous time-harmonic Maxwell equations posed on $\Omega$, which read as follows: For a given right-hand side $\solf$, find a vector field $\solu:\Omega\rightarrow\Co^3$ satisfying

\begin{subequations}\label{Maxwellorig}
\begin{equation}
		\curl\mu^{-1}\curl\solu-k^2\varepsilon\solu = \solf,
		\tag{\ref{Maxwellorig}}
\end{equation}
where $k\in\Co\setminus\{0\}$ is the wavenumber, $\mu^{-1}$ and $\varepsilon$ are complex-valued tensor fields which satisfy a coercivity condition and are piecewise $\CM{\ell+1}$-regular, but may be discontinuous across the interfaces $\interf_1,\ldots,\interf_r$.

Concerning boundary conditions on $\Gamma$, we allow for the following three choices:
\begin{itemize}
		\item Inhomogeneous impedance boundary conditions on $\Gamma$, which read as 
				\begin{equation}
						\left(\mu^{-1}\curl\solu\right)\times\soln-ik\zeta\solu_T = \solgi,
						\label{Maxwellimpedance}
				\end{equation}
				for a given tangent field $\solgi$, where $i:=\sqrt{-1}$ is the imaginary unit, $\soln$ is the outer unit normal to $\Gamma$, $\solu_T := \soln\times(\solu\times\soln)$ and $\zeta:\Gamma\rightarrow\Co^{3\times 3}$ is a given tensor field which is $\CM{\ell+1}$-regular, satisfies a coercivity condition as well as $(\zeta\solv)_T=\zeta\solv_T$ for all vector fields $\solv$ on $\Gamma$.

		\item Inhomogeneous natural boundary conditions on $\Gamma$, which require that
				\begin{equation}
						\left(\mu^{-1}\curl\solu\right)\times\soln = \solgn
						\label{Maxwellnatural}
				\end{equation}
				for a given tangent field $\solgn$.
				
		\item Homogeneous essential boundary conditions on $\Gamma$, which read as 
				\begin{equation}
				\solu_T = 0
				\label{Maxwellessential}
				\end{equation}
				on $\Gamma$, where, again, $\solu_T:=\soln\times(\solu\times\soln)$ denotes the tangential component of $\solu$.
\end{itemize}
\end{subequations}

We apply our findings concerning the regularity of vector fields with piecewise regular curl and divergence to \eqref{Maxwellorig} coupled by either \eqref{Maxwellimpedance}, \eqref{Maxwellnatural} or \eqref{Maxwellessential} and derive wavenumber-explicit regularity estimates for the corresponding solution $\solu$. 

Especially in the last decade, such wavenumber-explicit estimates have become increasingly important for the purposes of numerical analysis of Maxwell's equations (see e.g. \cite{MaxwellTomezyk,MaxwellImpedanceMelenk, MaxwellTransparentMelenk, MaxwellChaumont, MaxwellSpence}), hence our work as well as further research on this topic seem important.

\medskip 

So far, the regularity properties of Maxwell's equations was discussed in many papers, \cite{MaxwellTransparentMelenk, MaxwellImpedanceMelenk, MaxwellLu, MaxwellHiptmair, RegularityWeber, MaxwellChen, MaxwellCostabel}, to name only a few. However, many of these works either provide estimates that are either not wavenumber-explicit or they impose much stricter conditions on the coefficients and the geometry than we do in this work. 

The works that probably come closest to our work in terms of regularity estimates for Maxwell's equations are \cite{MaxwellChaumont} and \cite{MaxwellSpence}. In the former, however, the authors consider only homogeneous essential boundary conditions, real-valued tensors $\mu^{-1}$, $\varepsilon$ and divergence-free right-hand sides $\solf$. In the latter, the authors allow for complex-valued $\mu^{-1}$ and $\varepsilon$, but discuss only the case of homogeneous essential boundary conditions and solenoidal right-hand sides $\solf$. 

As of yet, there seems to be no paper dealing with wavenumber-explicit estimates for Maxwell's equations in the very general setting considered by us. Indeed, Theorem~\ref{Mainresult3} below provides regularity shift results for solutions $\solu$ of Maxwell's equations and meticulously tracks the wavenumber and the influence of $\solf$, $\diverg\solf$ and the boundary data in the corresponding regularity estimate. 

\medskip

The outline of this paper is as follows: In Section~\ref{notation} we describe the notation and assumptions considered throughout this work in greater detail and formulate our main results, namely Theorem~\ref{Mainresult1}, Theorem~\ref{Mainresult2} and Theorem~\ref{Mainresult3}. 

In Section~\ref{diffopsection} we recall the notions and some properties of surface differential operators before rewriting Maxwell's equations in a variational form. 

In Section~\ref{transmissionproblems} we develop the tools that we rely upon throughout Section~\ref{sec:regdecomp}. More precisely, Section~\ref{transmissionproblems} is about the regularity of certain Poisson transmission problems and the existence of Helmholtz decompositions for vector fields that are only piecewise regular. 

In Section~\ref{sec:regdecomp} we give an elegant proof for the fact that the (piecewise) Sobolev norm of order $\ell+1$ of a vector field $\solv$ can be controlled by the (piecewise) Sobolev norms of order $\ell$ of $\diverg\solv$ and $\curl\solv$ plus some boundary term, thus showing Theorem~\ref{Mainresult1} and Theorem~\ref{Mainresult2}.

Finally, in Section~\ref{finalsec}, we provide a proof for our third main result, namely Theorem~\ref{Mainresult3}.

\medskip 

Let us mention that this work is the first paper of an upcoming \revision{series} of research articles. The second paper \cite{MaxwellMyself2} in this \revision{series} builds on the subsequent Theorem~\ref{Mainresult1} and Theorem~\ref{Mainresult2} and proves wavenumber-explicit piecewise analyticity of solutions $\solu$ of \eqref{Maxwellorig}, provided that the geometry and all given coefficients and data are (piecewise) analytic. \revision{Subsequently}, the third work \cite{MaxwellMyself3} of this upcoming \revision{series} provides a wavenumber-explicit analysis of the $hp$-finite element method applied to Maxwell's equations \eqref{Maxwellorig} in the presence of piecewise analytic coefficients. In particular, if $k$ is the wavenumber, $h$ is the local mesh-width and $p$ is the local polynomial degree in the $hp$-FEM approximation, the third paper will show that the $hp$-FEM is quasi-optimal under the scale-resolution conditions a) that $|hk/p|$ is sufficiently small and b) that $p/\log |k|$ is sufficiently large.

%% file: notation.tex
\section{Notation and main results}\label{notation}

For any two vectors $\solw,\solz\in\Co^3$ with $\solw = (w_1, w_2, w_3)^T$ and $\solz=(z_1,z_2,z_3)^T$ we set $\solw\cdot\solz:=\sum_{i=1}^3w_iz_i$ and write $\SCP{\solw}{\solz}{}:=\solw\cdot\overline{\solz}$ for the scalar product between $\solw$ and $\solz$, where $\overline{\solz}:=(\overline{z_1},\overline{z_2},\overline{z_3})^T$ denotes the complex conjugate of $\solz$. 
Furthermore, the cross product between the vectors $\solw$ and $\solz$ is defined in the usual way as $\solw\times\solz:=(w_2z_3-w_3z_2, w_3z_1-w_1z_3, w_1z_2-w_2z_1)^T$.
  
\medskip 

  As usual, let $\SL(\Omega)$ denote the Lebesgue space of complex-valued square integrable functions, and define its vector-valued version $\VSL(\Omega):=(\SL(\Omega))^3$.
  For a bounded Lipschitz domain $\Omega\subseteq \R^3$ and $s\geq 0$, the space $\SHM{s}(\Omega)$ is the usual Sobolev space of order $s$, see \cite[Chapter 3]{BookMcLean}, and $\SHM{s}_0(\Omega)$ denotes the closure of $\CM{\infty}_0(\Omega)$ in $\SHM{s}(\Omega)$. \revision{Furthermore, in order to deal with vector fields, we define the vector-valued space $\VSHM{s}(\Omega):=(\SHM{s}(\Omega))^3$ and set}
  \begin{align*}
  \Hcurl := \{\solu\in\VSL(\Omega)\ |\ \curl \solu\in\VSL(\Omega)\}\quad {\rm and}\quad \Hdiv:=\{\solu\in\VSL(\Omega)\ |\ \diverg \solu\in\SL(\Omega)\}.
  \end{align*}

  \bigskip

  Finally, for a sufficiently smooth \revision{and orientable} surface $\Sigma$ and $s> 0$, let $\SHM{s}(\Sigma)$ be the fractional Sobolev space of index $s$ with dual space $\SHM{-s}(\Sigma)$, see \cite[Chapter 3]{BookMcLean}, and let $\VSHM{s}(\Sigma):=(\SHM{s}(\Sigma))^3$. Furthermore, \revision{since $\Sigma$ is supposed to be orientable we may choose a normal unit vector field $\soln$ on $\Sigma$ and define the space of square-integrable tangent fields by}
\begin{align*}
\VSL_T(\Sigma):=\{\solv\in\VSL(\Sigma)\ |\ \solv\cdot\soln = 0\}.
\end{align*}
\revision{For $s\geq 0$ we set}
\begin{align*}
\VSHM{s}_T(\Sigma):=\VSL_T(\Sigma)\cap\VSHM{s}(\Sigma),\ \ {\rm as\ well\ as}\ \ \VSHM{-s}_T(\Sigma):=(\VSHM{s}_T(\Sigma))'.
\end{align*}

  For a bounded Lipschitz domain $\Omega\subseteq\R^3$ the symbol $\overline{\Omega}$ denotes the closure of $\Omega$, and for $\ell\in\N\cup\{\infty\}$, a complex-valued function $v$ on $\Omega$ is in $\CM{\ell}(\overline{\Omega})$ if $v$ can be extended to a function $\widetilde{v}\in\CM{\ell}(\R^3)$.
  Similarly, vector fields $\solv$ and tensor fields $\nu$ on $\Omega$ are in $\VCM{\ell}(\overline{\Omega})$ if all of their respective component functions are in $\CM{\ell}(\overline{\Omega})$.

\subsection{$\CM{\ell}$-partitions and coercive $\VCMP{m}$-tensor fields}

Let $\Omega\subseteq\R^3$ be a bounded Lipschitz domain with boundary $\Gamma$. 
Throughout this work we suppose that $\Omega$ is partitioned into $\Gp_1,\ldots,\Gp_n$ subdomains such that the boundaries of these subdomains form a set of mutually disjoint closed\footnote{By {\it closed surface} we denote a compact surface without boundary.} surfaces inside of $\Omega$. 
The following definition makes this precise.

\begin{definition}\label{partitiondef}
		Let $\ell\in\N\cup\{\infty\}$. A $\CM{\ell}$-partition $\Gp$ is a tuple $\Gp = \geom$ which consists of domains $\Omega, \Gp_1,\ldots, \Gp_n$ satisfying 
		\begin{itemize}
				\item[(i)] \revision{The domains $\Omega, \Gp_1,\ldots,\Gp_n$ are bounded three-dimensional $\CM{\ell}$-domains and $\Omega$ is simply connected. Furthermore, the domains $\Gp_1,\ldots\Gp_n$ are mutually disjoint and satisfy $\overline{\Omega} = \overline{\Gp_1}\cup\ldots\cup\overline{\Gp_n}$.}
				\item[(ii)] The boundary $\Gamma:=\partial\Omega$ consists of $n'\geq 1$ simply connected components $\Gamma_1,\ldots, \Gamma_{n'}$. 
				\item[(iii)] There exist closed and simply connected $\CM{\ell}$-surfaces $\interf_1,\ldots, \interf_r$ such that $\Gamma, \interf_1,\ldots,\interf_r$ are mutually disjoint and such that 
\begin{align*}
		\Gamma\cup\bigcup_{j=1}^r\interf_j = \bigcup_{j=1}^n\partial\Gp_j.
\end{align*}
The surfaces $\interf_1,\ldots,\interf_r$ are called interface components and their union $\interf := \interf_1\cup\ldots\cup\interf_r$ is referred to as the subdomain interface.
\end{itemize}
\end{definition}

The point of introducing subdomains in the preceding fashion is to incorporate the location of (possible) discontinuities of piecewise regular vector fields into the geometry of the problem. 

We notice that requirement (iii) in Definition~\ref{partitiondef} implies that every connected component of $\partial\Gp_j$ either coincides with $\Gamma$ or with an interface component $\interf_i$. Moreover, for every $\interf_i$ there are precisely two subdomains $\Gp_j$ and $\Gp_h$ such that $\interf_i = \Gp_j\cap\Gp_h$. 
In particular, there may be no point in $\Omega$ where three or more subdomains meet.

\revision{
		\begin{remark}\label{nsc1}
				The assumption that $\Omega$ is simply connected makes some proofs in this work simpler (especially the proofs of Lemma~\ref{reg:helmholtze} and Lemma~\ref{reg:helmholtzdiv}), but it is not fundamental for our results. Indeed, the regularity results of Theorem~\ref{Mainresult1} and Theorem~\ref{Mainresult2} below hold true for a non-simply connected domain $\Omega$ as well with only minor modifications necessary. 
				The difference between the simply connected case considered in the following and the non-simply connected case is that in the latter the cohomology spaces of the de Rham complex have to be taken into account. It is known that on a bounded Lipschitz domain, these cohomology spaces are finite dimensional spaces of smooth functions or vector fields \cite{CurlInverse}. Hence, the existence of non-trivial cohomology spaces does not spoil the regularity results that we derive for the case of simply connected $\Omega,\Gp_1,\ldots,\Gp_n$, only the norm bounds have to be adapted. We come back to this issue in Remark~\ref{nsc2} and Remark~\ref{nsc3}.
		\end{remark}

}

\medskip

The subsequent two definitions clarity the notion of piecewise regular functions, vector- and tensor fields.

\begin{definition}
		Let $\Gp$ be a $\CM{\ell}$-partition and $m\in\N\cup\{\infty\}$. The space $\CMP{m}(\overline{\Omega})$ consists of all $v:\Omega\rightarrow\Co$ such that for $i=1,\ldots,n$ the restriction of $v$ to $\Gp_i$ is an element of $\CM{m}(\overline{\Gp_i})$. 

		Furthermore, a vector field $\solv:\Omega\rightarrow\Co^3$ is in $\VCMP{m}(\overline{\Omega})$ if its component functions are elements of $\CMP{m}(\overline{\Omega})$.
\end{definition}

\begin{definition}\label{assumptioncoefs}
		Let $\Gp$ be a $\CM{\ell}$-partition and let $m\in\N\cup\{\infty\}$. A tensor field $\nu:\Omega\rightarrow\Co^{3\times 3}$ is called a $\VCMP{m}$-tensor field, if its components are in $\CMP{m}(\overline{\Omega})$. 

		In addition, a $\VCMP{m}$-tensor field $\nu$ is called coercive, if there exist $\alpha\in\Co$ with $|\alpha|=1$ and $c>0$ such that 
		\begin{align*}
				\realpart\SCP{\alpha\nu\solz}{\solz}{}\geq c\norm{\solz}{}^2	
		\end{align*}
		for all $\solz\in\Co^3$ uniformly on $\Omega$.
\end{definition}

For the treatment of Maxwell's equations with impedance boundary conditions we will need the notion of $\VCM{m}$-tensor fields acting in the tangent plane to an orientable $\CM{\ell}$-surface $\Sigma$. For a vector field $\solv$ on $\Sigma$ we write $\solv_T := \soln\times(\solv\times\soln)$ for the tangential component of $\solv$, where $\soln$ is a unit normal to $\Sigma$.

\begin{definition}\label{assumptionimpedance}
		Let $\Sigma$ be an orientable $\CM{\ell}$-surface and suppose $m\leq \ell$. 
	A tensor field $\lambda:\Sigma\rightarrow\Co^{3\times 3}$ is called a $\VCM{m}$-tensor field acting in the tangent plane to $\Sigma$ if all components of $\lambda$ are in $\CM{m}(\Sigma)$ and there holds $(\lambda\solv)_T = \lambda\solv_T$ for all vector fields $\solv$ on $\Sigma$.

	A $\VCM{m}$-tensor field $\lambda$ acting in the tangent plane to $\Sigma$ is called coercive, if there exist $\alpha\in\Co$ with $|\alpha|=1$ and $c>0$ such that
		\begin{align*}
				\realpart\SCP{\alpha\lambda\solz}{\solz}{}\geq c\norm{\solz}{}^2	
		\end{align*}
		for all $\solz\in\Co^3$ uniformly on $\Sigma$.

\end{definition}

In the following we define piecewise Sobolev spaces, which play a crucial role in this work. To that end, let $\Gp$ be a $\CM{\ell}$-partition and assume $m\in\N_0$. Then, we define the spaces
\begin{align*}
		\PSHM{m}(\Gp):=\{u\in\SL(\Omega)\ \big|\ u\vert_{\Gp_i}\in\SHM{m}(\Gp_i)\ {\rm for}\ i=1,\ldots,n\} \quad{\rm and}\quad \PVSHM{m}(\Gp):=(\PSHM{m}(\Gp))^3,
\end{align*}
which are equipped with the norms
\begin{align*}
		\norm{u}{\PSHM{m}(\Gp)}^2:=\sum_{j=1}^n\norm{u}{\SHM{m}(\Gp_j)}^2\quad{\rm and} \quad \norm{\solu}{\PVSHM{m}(\Gp)}^2:=\sum_{i=1}^3\norm{\solu_i}{\PSHM{m}(\Gp)}^2,
\end{align*}
respectively, 
where $\solu_i$ denotes the $i$-th component of $\solu$.
Finally, for $m\in\N_0$ we introduce the spaces of vector fields with piecewise regular curl and divergence
\begin{align*}
\PVHcurl{m}:=\{\solv\in\VSL(\Omega)\ \big|\ \curl\solv\vert_{\Gp_i}\in\VSHM{m}(\Gp_i)\ {\rm for}\ i=1,\ldots,n\}
\end{align*}
and
\begin{align*}
\PVHdiv{m}:=\{\solv\in\VSL(\Omega)\ \big|\ \diverg\solv\vert_{\Gp_i}\in\SHM{m}(\Gp_i)\ {\rm for}\ i=1,\ldots,n\}.
\end{align*}
\revision{Notice that we defined $\PVHcurl{m}$ and $\PVHdiv{m}$ in a rather unconventional way. Canonically, one would rather write $\solv\in\PVHcurl{m}$ if $\solv\in\PVSHM{m}(\Gp)$ and $\curl\solv\in\PVSHM{m}(\Gp)$ and analogously for $\PVHdiv{m}$. We deliberately chose this unorthodox way in order to be able to talk about (piecewise) regularity of $\curl\solu$ and $\diverg\solu$ independently of the (piecewise) regularity of $\solu$.}

\subsection{Main results}

The first main result of our work deals with the piecewise regularity of a vector field $\solu$ with piecewise regular curl and divergence, as well as possibly inhomogeneous normal trace on $\Gamma$. It shows that for such a vector field there exists an $\VSL(\Omega)$-orthogonal decomposition into a gradient field and a vector field which is solely controlled by the curl of $\solu$.

We highlight that compared to \cite[Theorem 2.2]{RegularityWeber} or \cite[Theorem 9.1]{MaxwellSpence}, our result allows for inhomogeneous normal traces on $\Gamma$ and is also slightly sharper, since our estimates involve only the curl, divergence and normal trace of $\solu$, but are independent of $\norm{\solu}{\VSL(\Omega)}$. \revision{This independence of $\norm{\solu}{\VSL(\Omega)}$ is due to our assumptions of $\Omega$ being simply connected and of $\Gamma$ and $\interf$ consisting of simply connected components.}

\begin{theorem}\label{Mainresult1}
		Let $\ell\in\N_0$, suppose that $\Gp$ is a $\CM{\ell+2}$-partition and let $\nu$ be a coercive $\VCMP{\ell+1}$-tensor field in the sense of Definition~\ref{assumptioncoefs}.
Under these assumptions, suppose that $\solu\in\Hcurl$ satisfies $\solu\in\PVHcurl{\ell}$ and $\nu\solu\in\Hdiv\cap\PVHdiv{\ell}$, as well as $\nu\solu\cdot\soln = h$ on $\Gamma$ for a function $h\in\SHM{\ell+1/2}(\Gamma)$. 

		 Then, $\solu\in\PVSHM{\ell+1}(\Gp)$, and there exists a decomposition $\solu = \solv+\nabla\phi$ with 
		 \revision{
				 \begin{align}\label{res1estimates}
						 \begin{split}
						  \norm{\solv}{\PVSHM{\ell+1}(\Gp)} &\leq C\norm{\curl\solu}{\PVSHM{\ell}(\Gp)}, \\
						 \norm{\phi}{\SHM{1}(\Omega)}+\norm{\phi}{\PSHM{\ell+2}(\Omega)}&\leq C\left(\norm{\diverg\nu\solu}{\PVSHM{\ell}(\Gp)}+\norm{h}{\SHM{\ell+1/2}(\Gamma)}\right),
						 \end{split}
				 \end{align}
		 and} $\SCP{\nu\solv}{\nabla\xi}{\VSL(\Omega)}=0$ for all $\xi\in\SHM{1}(\Omega)$. Moreover, the constant $C>0$ depends only on $\Gp$, $\nu$, and $\ell$.


\end{theorem}

Our second main result can be seen as a dual statement of Theorem~\ref{Mainresult1}; it deals with possibly inhomogeneous tangential traces of $\solu$. Again, it states the existence of a Helmholtz-type decomposition, this time into a curl-field and something which can be controlled by the divergence of $\solu$. 

Again, we highlight that our estimates involve only the curl and divergence of $\solu$ and the tangential component $\solu_T:=\soln\times(\solu\times\soln)$ or the tangential trace $\solu_t:=\solu\times\soln$, but are independent of the $\VSL(\Omega)$-norm of $\solu$.

\begin{theorem}\label{Mainresult2}
		Let $\ell\in\N_0$, suppose that $\Gp$ is a $\CM{\ell+2}$-partition and let $\nu$ be a coercive $\VCMP{\ell+1}$-tensor field in the sense of Definition~\ref{assumptioncoefs}.
		Under these assumptions, suppose that $\solu\in\Hcurl$ satisfies $\solu\in\PVHcurl{\ell}$ and $\nu\solu\in\Hdiv\cap\PVHdiv{\ell}$, as well as either $\solu_T = \solg_T$ on $\Gamma$ or $\solu_t=\solg_T$ on $\Gamma$ for a tangent field $\solg_T\in\VSHM{\ell+1/2}_T(\Gamma)$. 

	Then, $\solu\in\PVSHM{\ell+1}(\Gp)$, and there exists a decomposition $\nu\solu = \nu\solv+\curl\solz$ with 
	\revision{
			\begin{align}\label{res2estimates}
					\begin{split}
						\norm{\solv}{\PVSHM{\ell+1}(\Gp)} &\leq C\norm{\diverg\nu\solu}{\PVSHM{\ell}(\Gp)},\\
						\norm{\solz}{\VSHM{1}(\Omega)}+\norm{\curl\solz}{\PVSHM{\ell+1}(\Omega)}&\leq C\left(\norm{\curl\solu}{\PVSHM{\ell}(\Gp)}+\norm{\solg_T}{\VSHM{\ell+1/2}_T(\Gamma)}\right),
					\end{split}
							\end{align}
	and} $\SCP{\solv}{\curl\solw}{\VSL(\Omega)}=0$ for all $\solw\in\Hcurl$.
	Moreover, the constant $C>0$ depends only on $\Gp$, $\nu$, and $\ell$.


\end{theorem}

\revision{

		\begin{remark}
				In the case that the $\CM{\ell+2}$-partition $\Gp$ consists only of one subdomain (that is, there are no surfaces of discontinuity of $\nu$), the statements of Theorem~\ref{Mainresult1} and Theorem~\ref{Mainresult2} hold with the broken Sobolev norms replaced by their "usual" counterparts: If $\nu$ satisfies $\nu\in\VCM{\ell+1}(\overline{\Omega})$ and $\solu\in\VHcurl{\ell}\cap\VHdiv{\ell}$, then the statements of Theorem~\ref{Mainresult1} and Theorem~\ref{Mainresult2} hold true with $\PVSHM{\ell+1}(\Gp)$ and $\PVSHM{\ell}(\Gp)$ replaced by $\VSHM{\ell+1}(\Omega)$ and $\VSHM{\ell}(\Omega)$ in all instances.
		\end{remark}
		
}

\revision{
		\begin{remark}\label{nsc2}
				In the case that $\Omega$ is not simply connected or if $\Gamma$ or the subdomain interface $\interf$ have non-simply connected components, Theorem~\ref{Mainresult1} and Theorem~\ref{Mainresult2} still remain valid once we change \eqref{res1estimates} and \eqref{res2estimates} to
\begin{align*}
		\norm{\solv}{\PVSHM{\ell+1}(\Gp)} &\leq C\left(\norm{\solu}{\VSL(\Omega)}+\norm{\curl\solu}{\PVSHM{\ell}(\Gp)}\right), \\
		\norm{\phi}{\SHM{1}(\Omega)}+\norm{\phi}{\PSHM{\ell+2}(\Omega)}&\leq C\left(\norm{\solu}{\VSL(\Omega)}+\norm{\diverg\nu\solu}{\PVSHM{\ell}(\Gp)}+\norm{h}{\SHM{\ell+1/2}(\Gamma)}\right),
				 \end{align*}
		and
	\begin{align*}
			\norm{\solv}{\PVSHM{\ell+1}(\Gp)} &\leq C\left(\norm{\solu}{\VSL(\Omega)}+\norm{\diverg\nu\solu}{\PVSHM{\ell}(\Gp)}\right),\\
			\norm{\solz}{\VSHM{1}(\Omega)}+\norm{\curl\solz}{\PVSHM{\ell+1}(\Omega)}&\leq C\left(\norm{\solu}{\VSL(\Omega)}+\norm{\curl\solu}{\PVSHM{\ell}(\Gp)}+\norm{\solg_T}{\VSHM{\ell+1/2}_T(\Gamma)}\right),
							\end{align*}
							respectively, cf. Remark~\ref{nsc1}.	
		\end{remark}
}

The third and final main result of this paper follows from applying the previous two theorems to the time-harmonic Maxwell equations \eqref{Maxwellorig}. It provides piecewise Sobolev regularity of a weak solution $\solu$ together with wavenumber-explicit estimates.

The influence of the wavenumber $k$ on the regularity of $\solu$ depends in a more complicated way on different quantities related to the given data $\solf$, $\solgn$ and $\solgi$. To ease the presentation, we abbreviate
\begin{align*}
		\solF_{m,k} := |k|^{-m}\left(\norm{\solf}{\PVSHM{m}(\Gp)}+|k|^{-1}\norm{\diverg\solf}{\PSHM{m}(\Gp)}\right)
\end{align*}
for $m\in\N_0$, and, depending on the imposed boundary condition,
\begin{alignat*}{2}
		\solG_{m,k} &:= |k|^{-m}\left(\norm{\diverg\solf}{\PSHM{m}(\Gp)}+\norm{\solgi}{\revision{\VSHM{m+1/2}_T(\Gamma)}}\right)\quad && {\rm in\ the\ case\ of\ \eqref{Maxwellimpedance}}, \\ 
		\solG_{m,k} &:= |k|^{-m}\left(\norm{\solgn}{\VSHM{m+1/2}_T(\Gamma)}+|k|^{-1}\norm{\revision{\solf\cdot\soln}-\diverg_{\Gamma}\solgn}{\SHM{m+1/2}(\Gamma)}\right)\quad && {\rm in\ the\ case\ of\ \eqref{Maxwellnatural}},\\
		\solG_{m,k} &:= 0\quad  && {\rm in\ the\ case\ of\ \eqref{Maxwellessential}},
\end{alignat*}
where $\diverg_{\Gamma}$ denotes the surface divergence on $\Gamma$, see e.g. \cite{BookMonk, BookNedelec} or Section~\ref{diffopsection} below. 

\medskip

With these definitions, and with $\diverg_{\Gamma}$ being the aforementioned surface divergence, the third and last main result of this paper reads as follows:
\begin{theorem}\label{Mainresult3}
		Let $\ell\in\N_0$, suppose that $\Gp$ is a $\CM{\ell+2}$-partition and let $\mu^{-1}$ and $\varepsilon$ be coercive $\VCMP{\ell+1}$-tensor fields in the sense of Definition~\ref{assumptioncoefs}. Moreover, if we impose impedance boundary conditions \eqref{Maxwellimpedance}, we suppose that $\zeta$ is a coercive $\VCM{\ell+1}$-tensor field acting in the tangent plane to $\Gamma$, see Definition~\ref{assumptionimpedance}.

		\medskip 

		Under these assumptions, let $\solu$ be a weak solution of \eqref{Maxwellorig} in the sense of Section~\ref{secweakform} below, and suppose that the right-hand side $\solf\in\Hdiv$ satisfies $\solf\in\PVSHM{\ell}(\Gp)\cap\PVHdiv{\ell}$. Depending on the imposed boundary conditions, we assume 

		\begin{itemize}
				\item $\solgi\in\VSHM{\ell+1/2}_T(\Gamma)$ in the case of impedance boundary conditions \eqref{Maxwellimpedance}.
				\item $\solgn\in\VSHM{\ell+1/2}_T(\Gamma)$ and $\revision{\solf\cdot\soln-}\diverg_{\Gamma}\solgn\in\SHM{\ell+1/2}(\Gamma)$ in the case of natural boundary conditions \revision{\eqref{Maxwellnatural}}. 
		\end{itemize}

		Then, \revision{for all $m\in\{0,\ldots,\ell\}$} there holds $\solu\in\PVSHM{\revision{m}+1}(\Gp)\cap\PVHcurl{\revision{m}+1}$ with the estimate
	\begin{align}\label{finiteregularity}
			\begin{split}
					\revision{|k|^{-m-1}}\norm{\solu}{\PVSHM{\revision{m}+1}(\Gp)}+\revision{|k|^{-m-2}}\norm{\curl\solu}{\PVSHM{\revision{m}+1}(\Gp)}\leq \revision{C_m} &\left(\norm{\solu}{\VSL(\Omega)}+|k|^{\revision{-1}}\norm{\curl\solu}{\VSL(\Omega)}\right) \\
																																										  &\hskip 3cm +\revision{C_m}|k|^{\revision{-2}}\sum_{j=0}^{\revision{m}}\left[\solF_{j,k}+\solG_{j,k}\right],
	\end{split}
	\end{align}
	where the constant $\revision{C_m}>0$ depends only on $\Gp$, $\revision{m}$, $\mu^{-1}$, $\varepsilon$, and, if necessary, $\zeta$. 

%
	
\end{theorem}

\revision{
		\begin{remark}\label{nsc3}
				If $\Omega$ is not simply connected or if $\Gamma$ or the subdomain interface $\interf$ have non-simply connected components, Theorem~\ref{Mainresult3} still remains valid once we change \eqref{finiteregularity} to
\begin{align*}
			\begin{split}
					\revision{|k|^{-m-1}}\norm{\solu}{\PVSHM{\revision{m}+1}(\Gp)}+\revision{|k|^{-m-2}}\norm{\curl\solu}{\PVSHM{\revision{m}+1}(\Gp)}\leq \revision{C_m}& \left(1+\sum_{j=0}^m|k|^{-j-1}\right)\left(\norm{\solu}{\VSL(\Omega)}+|k|^{\revision{-1}}\norm{\curl\solu}{\VSL(\Omega)}\right) \\
																																										  &\hskip 3cm +C_m|k|^{\revision{-2}}\sum_{j=0}^{\revision{m}}\left[\solF_{j,k}+\solG_{j,k}\right].
	\end{split}
	\end{align*}
	For $|k|\geq k_0>0$ this reduces again to \eqref{finiteregularity} for a modified constant $C_{m}$ that depends also on $k_0$.
		\end{remark}
}

%
%
%
Let us \revision{highlight} that in general there are $k\in\Co\setminus\{0\}$, right-hand sides $\solf$ and, if necessary, boundary data $\solgi$ or $\solgn$ such that no weak solution $\solu$ of \eqref{Maxwellorig} exists. However, if there is a (possible not unique) weak solution $\solu$, then Theorem~\ref{Mainresult3} asserts piecewise Sobolev regularity of any weak solution of \eqref{Maxwellorig}.

\medskip

In addition, as a consequence of Theorem~\ref{Mainresult3} and the Sobolev embedding theorem we have the following corollary:

\begin{corollary}
		Suppose that $\Gp$ is a $\CM{\infty}$-geometry, that $\mu^{-1}$ and $\varepsilon$ are coercive $\VCMP{\infty}$-tensor fields in the sense of Definition~\ref{assumptioncoefs} and assume that $\solf\in\Hdiv$. Moreover,
		\begin{itemize}
				\item in the case of impedance boundary conditions \eqref{Maxwellimpedance}, let $\zeta$ be a coercive $\VCM{\infty}$-tensor field acting in the tangent plane to $\Gamma$ and let $\solgi$ be a smooth tangent field on $\Gamma$, 
				\item in the case of natural boundary conditions \eqref{Maxwellimpedance}, let $\solgn$ be a smooth tangent field to $\Gamma$. 
\end{itemize}
Under these assumptions, consider a weak solution $\solu$ of \eqref{Maxwellorig} in the sense of Section~\ref{secweakform} below. If $\solf$ is piecewise smooth, i.e., $\solf\in\VCMP{\infty}(\overline{\Omega})$, then \eqref{finiteregularity} holds for all $\ell\in\N_0$. Consequently, $\solu$ is piecewise smooth, i.e., there holds $\solu\in\VCMP{\infty}(\overline{\Omega})$.
\end{corollary}

%% file: diffops.tex
\section{Differential operators on surfaces and traces of $\Hcurl$}\label{diffopsection}

In this section we recollect some properties of surface differential operators and recall the interplay between these operators and the canonical trace operators on $\Hcurl$. 

For most of this section we assume that $\Sigma$ is a closed \revision{(i.e., compact and without boundary)} and orientable $\CM{2}$-surface consisting of $n'\geq 1$ simply connected components\footnote{You can think of $\Sigma$ being a simply connected component $\Gamma_j$ of $\Gamma$ or an interface component $\interf_j$.}, and $\soln:\Sigma\rightarrow\mathbb{S}_2$ denotes a unit vector field normal to $\Sigma$. \revision{Let us mention that for many statements of this section, the assumption of $\Sigma$ being $\CM{2}$ is stronger than necessary. Indeed, many results can be extended to the case of $\Sigma$ being only Lipschitz. However, for Lipschitz surfaces some statements become more technical as one has to be more careful in the analysis of the subsequently defined surface differential operators. Hence, in order to make things easier we consider only surfaces that are at least $\CM{2}$.}

We briefly recall some definitions and results concerning surface differential operators from \cite{BookMonk, BookNedelec}. Let $\Sigma_{\tau}$ be a sufficiently small tubular neighborhood around $\Sigma$. Following the notation from \cite{MaxwellImpedanceMelenk,BookNedelec}, the constant extensions (in normal direction) of a sufficiently smooth scalar function $u$ on $\Sigma$ is denoted by $u^*$; the surface gradient $\nabla_{\Sigma}$ and the tangential curl operator  $\overrightarrow{\curl_{\Sigma}}$ are then defined as 
\begin{align*}
\nabla_{\Sigma}:= (\nabla u^*)\vert_{\Sigma}\quad {\rm and}\ \ \overrightarrow{\curl_{\Sigma}}u := \nabla_{\Sigma}u\times \soln.
\end{align*} 
Note that $-\nabla_{\Sigma}$ and $\overrightarrow{\curl_{\Sigma}}$ map scalar functions to tangent fields. Thus, their adjoint operators $\diverg_{\Sigma}$ and $\curl_{\Sigma}$ map tangent fields to scalar functions and satisfy 
\begin{align*}
	\SCP{\nabla_{\Sigma} u}{\solv}{\VSL(\Sigma)} = -\SCP{u}{\diverg_{\Sigma}\solv}{\SL(\Sigma)}\quad {\rm and}\quad \SCP{\overrightarrow{\curl_{\Sigma}}u}{\solv}{\VSL(\Sigma)}= \SCP{u}{\curl_{\Sigma}\solv}{\SL(\Sigma)}
\end{align*}
for all sufficiently smooth scalar functions $u$ and sufficiently smooth tangent fields $\solv$. 

\begin{remark}
		For simplicity we chose to introduce $\diverg_{\Sigma}$ and $\curl_{\Sigma}$ as adjoint operators to $-\nabla_{\Sigma}$ and $\overrightarrow{\curl_{\Sigma}}$, respectively. For a more rigorous approach we refer e.g. to \cite[Section 2.5.6]{BookNedelec}. Moreover, the definition of the above surface differential operators can be extended to Lipschitz surfaces, \cite{TracesHcurl}.
\end{remark}

The following lemma shows that the surface gradient $\nabla_{\Sigma}$ and the tangential curl operator $\overrightarrow{\curl_{\Sigma}}$ are connected to traces of volume gradients. Similarly, the surface divergence $\diverg_{\Sigma}$ and the surface curl operator $\curl_{\Sigma}$ are connected to traces of volume curls.
\begin{lemma}\label{setting:sglemma}
		Let $\Omega\subseteq\R^3$ be a bounded $\CM{2}$-domain with boundary $\Sigma$ and assume that $\Sigma$ consists of $n'\geq 1$ simply connected components. Let $\soln$ be the outer unit normal field to $\Sigma$. 
		Then, for any sufficiently smooth function $\psi$ there holds
	\begin{align*}
		\nabla_{\Sigma}(\psi\vert_{\Sigma}) = \soln\times(\nabla\psi\vert_{\Sigma}\times\soln) \quad{\rm and}\quad \overrightarrow{\curl_{\Sigma}}(\psi\vert_{\Sigma}) = \nabla\psi\vert_{\Sigma}\times\soln.
	\end{align*}
	
	Moreover, for any sufficiently smooth vector field $\solv$ there holds
	\begin{align}\label{surfacecurlconnect}
		\diverg_{\Sigma}\solv_t = -\curl_{\Sigma}\solv_T = (\curl\solv)\vert_{\Sigma}\cdot\soln,
	\end{align}
	where $\solv_t:=\solv\vert_{\Sigma}\times\soln$ and $\solv_T:=\soln\times(\solv\vert_{\Sigma}\times\soln)$.
\end{lemma}
\begin{fatproof}
	The statement concerning the surface gradient is \cite[(2.26)]{MaxwellImpedanceMelenk}. We note that this already implies $\overrightarrow{\curl_{\Sigma}}(\psi\vert_{\Sigma}) = \nabla\psi\vert_{\Sigma}\times\soln$.
	
	According to \cite[Theorem 3.31]{BookMonk} there holds $\SCP{\curl\solv}{\nabla\xi}{\VSL(\Omega)} = \SCP{\solv_t}{\nabla_{\Sigma}\xi}{\VSL(\Sigma)}$ for all sufficiently smooth functions $\xi$. Thus, partial integration and the definition of $\diverg_{\Sigma}$ yield $\SCP{\curl\solv\cdot\soln}{\xi}{\SL(\Sigma)} = \SCP{\diverg_{\Sigma}\solv_t}{\xi}{\VSL(\Sigma)}$, which shows $\diverg_{\Sigma}\solv_t = (\curl\solv)\vert_{\Sigma}\cdot\soln$.
	
	Furthermore, we observe that $\SCP{\curl\solv}{\nabla\xi}{\VSL(\Omega)} =\SCP{\solv_t}{\nabla_{\Sigma}\xi}{\VSL(\Sigma)} = \SCP{\solv_t\times\soln}{\overrightarrow{\curl_{\Sigma}}\xi}{\VSL(\Sigma)}$. Similarly as above, partial integration and the definition of $\curl_{\Sigma}$ yield $(\curl\solv)\vert_{\Sigma}\cdot\soln=-\curl_{\Sigma}\solv_T $, which concludes the proof.
\end{fatproof}

The following mapping property is proved in \cite[Proposition 3.6]{TracesHcurl} for bounded Lipschitz domains and $\ell=0$. For more regular domains we can extend this to the following result:
\begin{proposition}\label{setting:regularitysurfaceops}
		\revision{For any $\ell\in\N$ and $s\in [-\ell,\ell+1]$ and any closed and orientable $\CM{\ell+1}$-surface} $\Sigma$ consisting of $n'\geq 1$ simply connected components, the surface differential operators $\nabla_{\Sigma}$ and $\overrightarrow{\curl_{\Sigma}}$ extend to bounded linear operators
$\nabla_{\Sigma}:\SHM{\revision{s}}(\Sigma)\rightarrow \VSHM{\revision{s-1}}_T(\Sigma)$ and $\overrightarrow{\curl_{\Sigma}}:\SHM{\revision{s}}(\Sigma)\rightarrow \VSHM{\revision{s-1}}_T(\Sigma)$. Moreover, the operators $\diverg_{\Sigma}$ and $\curl_{\Sigma}$ extend to bounded linear operators $\diverg_{\Sigma}:\VSHM{\revision{s}}_T(\Sigma)\rightarrow\SHM{\revision{s-1}}(\Sigma)$ and $\curl_{\Sigma}:\VSHM{\revision{s}}_T(\Sigma)\rightarrow\SHM{\revision{s-1}}(\Sigma)$.
\end{proposition}

\revision{\begin{remark}
				In the case $\ell=0$ and if $\Sigma$ is merely Lipschitz, the above statement is more involved, as the ranges of the surface differential operators cannot be characterized as easily by the "typical" Sobolev spaces on $\Sigma$. For surfaces that are at least $\CM{2}$ things become much easier, since the associated normal vector fields are then at least $\CM{1}$ which leads to the surface differential operators behaving very nicely. Nevertheless we stress that the regularity assumptions on $\Sigma$ from Proposition~\ref{setting:regularitysurfaceops} are not optimal and can probably be relaxed by a more thorough analysis. 
	
				Similarly, the assumption that $\Sigma$ is $\CM{2}$ from Lemma~\ref{setting:sglemma} and Proposition~\ref{setting:Hodgedecomp} below is probably stronger than necessary. However, for our purposes the assumption that $\Sigma$ is $\CM{2}$ is sufficient, hence we do not aim to extend these statements to less regular surfaces.
\end{remark}}

The next proposition is a fundamental result from \cite[Section 3 and Section 5]{TracesHcurl}. It implies an exact sequence property of surface differential operators and asserts the existence of a Hodge decomposition for tangent fields. In \cite{TracesHcurl} it is stated for the more general case of boundaries of Lipschitz domains, however, since in this work we consider more regular domains, we formulate it for $\CM{2}$-surfaces.
\begin{proposition}\label{setting:Hodgedecomp}
		Let $\Sigma$ be a closed, simply connected and orientable $\CM{2}$-surface. Then there holds
	\begin{align*}
	\operatorname{ker}(\curl_{\Sigma})\cap\VSL_T(\Sigma) = \nabla_{\Sigma}\SHM{1}(\Sigma)\quad &{\rm and}\ \ \operatorname{ker}(\diverg_{\Sigma})\cap\VSL_T(\Sigma) = \overrightarrow{\curl_{\Sigma}}\SHM{1}(\Sigma),\\ 
	\operatorname{ker}(\curl_{\Sigma})\cap\VSHM{-1/2}_T(\Sigma) = \nabla_{\Sigma}\SHM{1/2}(\Sigma)\quad &{\rm and}\ \ \operatorname{ker}(\diverg_{\Sigma})\cap\VSHM{-1/2}_T(\Sigma) = \overrightarrow{\curl_{\Sigma}}\SHM{1/2}(\Sigma),
	\end{align*}
	as well as 
	\begin{align}\label{setting:kersurfacegrad}
	\operatorname{ker}(\nabla_{\Sigma})\cap \SHM{1/2}(\Sigma) = \operatorname{ker}(\overrightarrow{\curl_{\Sigma}})\cap \SHM{1/2}(\Sigma) = \R.
	\end{align}
	Moreover, there holds $\VSL_T(\Sigma) = \nabla_{\Sigma}\SHM{1}(\Sigma)\oplus\overrightarrow{\curl_{\Sigma}}\SHM{1}(\Sigma)$ and this decomposition is $\VSL_T$-orthogonal, that is, $\nabla_{\Sigma}\SHM{1}(\Sigma)\perp\overrightarrow{\curl_{\Sigma}}\SHM{1}(\Sigma)$ in $\VSL_T(\Sigma)$.
	
\end{proposition}

\begin{remark}
	Note that in Proposition~\ref{setting:Hodgedecomp} it is assumed that $\Sigma$ is simply connected. If $\Sigma$ consists of $n'\geq 1$ simply connected components, Proposition~\ref{setting:Hodgedecomp} stays valid once we change \eqref{setting:kersurfacegrad} to
	\begin{align*}
			\operatorname{ker}(\nabla_{\Sigma})\cap \SHM{1/2}(\Sigma) = \operatorname{ker}(\overrightarrow{\curl_{\Sigma}})\cap \SHM{1/2}(\Sigma) = \mathcal{M},
	\end{align*}
	where $\mathcal{M}$ is the $n'$-dimensional space spanned by the functions that are equal to one on a single simply connected component and zero on the others.
\end{remark}

\medskip 

The subsequent lemma provides a shift result for tangent fields. It turns out to be very useful for our purposes.
\begin{lemma}\label{shiftresultsurface}
		Let $\ell\in\N_0$ and let $\Sigma$ be a closed and orientable $\CM{\ell+2}$-surface consisting of $n'\geq 1$ simply connected components. Let $\lambda$ be a coercive $\CM{\ell+1}$-tensor field acting in the tangent plane to $\Sigma$, see Definition~\ref{assumptionimpedance}. 

	For any tangent field $\solv\in\VSL_T(\Sigma)$ satisfying $\diverg_{\Sigma}\solv\in\SHM{\ell-1/2}(\Sigma)$ and $\curl_{\Sigma}(\lambda\solv)\in\SHM{\ell-1/2}(\Sigma)$, there holds $\solv\in\VSHM{\ell+1/2}_T(\Sigma)$ with
	\begin{align*}
		\norm{\solv}{\VSHM{\ell+1/2}_T(\Sigma)}\leq C\left(\norm{\diverg_{\Sigma}\solv}{\SHM{\ell-1/2}(\Sigma)}+\norm{\curl_{\Sigma}(\lambda\solv)}{\SHM{\ell-1/2}(\Sigma)}\right),
	\end{align*}
	where the constant $C>0$ depends only on $\lambda$, $\ell$  and $\Sigma$. 

	Similarly, when $\solv\in\VSL_T(\Sigma)$ satisfies $\curl_{\Sigma}\solv\in\SHM{\ell-1/2}(\Sigma)$ and $\diverg_{\Sigma}(\lambda\solv)\in\SHM{\ell-1/2}(\Sigma)$, there holds 
$\solv\in\VSHM{\ell+1/2}_T(\Sigma)$ with
\begin{align*}
		\norm{\solv}{\VSHM{\ell+1/2}_T(\Sigma)}\leq C\left(\norm{\diverg_{\Sigma}(\lambda\solv)}{\SHM{\ell-1/2}(\Sigma)}+\norm{\curl_{\Sigma}\solv}{\SHM{\ell-1/2}(\Sigma)}\right),
\end{align*}
where, again, the constant $C>0$ depends only on $\lambda$, $\ell$, and $\Sigma$.
\end{lemma}
\begin{fatproof}
	We only consider the case $\diverg_{\Sigma}\solv\in\SHM{\ell-1/2}(\Sigma)$ and $\curl_{\Sigma}(\lambda\solv)\in\SHM{\ell-1/2}(\Sigma)$, the other case follows analogously. We employ Proposition \ref{setting:Hodgedecomp} to decompose $\solv$ into $\solv = \nabla_{\Sigma}\psi+\overrightarrow{\curl_{\Sigma}}\xi$. Due to $\diverg_{\Sigma}\overrightarrow{\curl_{\Sigma}}=0$ we have $\diverg_{\Sigma}\nabla_{\Sigma}\psi = \diverg_{\Sigma}\solv\in\SHM{\ell-1/2}(\Sigma)$, and according to elliptic regularity theory this implies $\psi\in\SHM{\ell+3/2}(\Sigma)$ with $\norm{\psi}{\SHM{\ell+3/2}(\Sigma)}\leq C\norm{\diverg_{\Sigma}\solv}{\SHM{\ell-1/2}(\Sigma)}$.

	In addition, there holds $\curl_{\Sigma}(\lambda\overrightarrow{\curl_{\Sigma}}\xi) = \curl_{\Sigma}(\lambda\solv)-\curl_{\Sigma}(\lambda\nabla_{\Sigma}\psi)\in\SHM{\ell-1/2}(\Sigma)$. We note that $\curl_{\Sigma}\lambda\overrightarrow{\curl_{\Sigma}}$ is an elliptic operator, therefore elliptic regularity theory guarantees 
	\begin{align*}
		\norm{\xi}{\SHM{\ell+3/2}(\Sigma)}\leq C\left(\norm{\curl_{\Sigma}(\lambda\solv)}{\SHM{\ell-1/2}(\Sigma)}+\norm{\diverg_{\Sigma}\solv}{\SHM{\ell-1/2}(\Sigma)}\right).
	\end{align*}
	This concludes the proof.
\end{fatproof}

\subsection{Trace spaces of $\Hcurl$}

We return to our original setting and assume that $\Omega$ is a bounded $\CM{\ell+2}$-domain with boundary $\Gamma$ for some $\ell\in\N_0$. We note that Lemma~\ref{setting:sglemma} suggests a connection between $\diverg_{\Gamma}$, $\curl_{\Gamma}$ and traces of volume curls. In fact, these two operators characterize the range of the canonical trace operators on $\Hcurl$. 
To make this precise, we follow \cite{BookMonk} and define the auxiliary spaces

\begin{align*}
		\Hcurlgamma{\Gamma}:=\left\{\solu\in\VSHM{-1/2}_T(\Gamma)\ \big|\ \curl_{\Gamma}\solu\in\SHM{-1/2}(\Gamma)\right\}\quad {\rm with}\quad \norm{\solu}{\Hcurlgamma{\Gamma}}^2:=\norm{\solu}{\VSHM{-1/2}_T(\Gamma)}^2+\norm{\curl_{\Gamma}\solu}{\SHM{-1/2}(\Gamma)}^2
\end{align*}
and 
\begin{align*}
		\Hdivgamma{\Gamma}:=\left\{\solu\in\VSHM{-1/2}_T(\Gamma)\ \big|\ \diverg_{\Gamma}\solu\in\SHM{-1/2}(\Gamma)\right\}\quad {\rm with}\quad \norm{\solu}{\Hdivgamma{\Gamma}}^2:=\norm{\solu}{\VSHM{-1/2}_T(\Gamma)}^2+\norm{\diverg_{\Gamma}\solu}{\SHM{-1/2}(\Gamma)}^2.
\end{align*}

There holds the following trace result, see \cite[Theorem 3.29]{BookMonk} or \cite[Theorem 5.4.2]{BookNedelec}.

\begin{proposition}\label{traceprop}
		Let $\Omega\subseteq\R^3$ be a bounded Lipschitz domain with boundary $\Gamma$ and outer normal vector $\soln$. 
		We consider the maps $\Pi_{T}$ and $\Pi_{t}$, which for smooth vector fields $\solv$ on $\overline{\Omega}$ are defined as 
	\begin{align*}
			\Pi_{T}\solv := \soln\times(\solv\times\soln)\quad {\rm and}\quad \Pi_{t}\solv := \solv\times\soln.
	\end{align*}
	These maps extend to bounded  and surjective operators $\Pi_{T}:\Hcurl\rightarrow\Hcurlgamma{\Gamma}$ and $\Pi_{t}:\Hcurl\rightarrow\Hdivgamma{\Gamma}$. In addition, there exist bounded lifting operators $\liftcurl:\Hcurlgamma{\Gamma}\rightarrow\Hcurl$ and $\liftdiv:\Hdivgamma{\Gamma}\rightarrow\Hcurl$. 
\end{proposition}

\begin{remark}
In order to shorten notation, we abbreviate $\Pi_T\solv$ by $\solv_T$ and $\Pi_t\solv$ by $\solv_t$, respectively. 
\end{remark}

As a direct consequence of Proposition~\ref{traceprop} we have that the spaces $\Hdivgamma{\Gamma}$ and $\Hcurlgamma{\Gamma}$ are dual to each other: Indeed, for $\solv\in\Hdivgamma{\Gamma}$ and $\solw\in\Hcurlgamma{\Gamma}$ we may define a duality pairing $\dualpair{\solv}{\solw}{\Zdiv}{\Ycurl}$ by
\begin{align*}
\dualpair{\solv}{\solw}{\Zdiv}{\Ycurl}:=\SCP{\curl\liftdiv\solv}{\liftcurl\solw}{\VSL(\Omega)}-\SCP{\liftdiv\solv}{\curl\liftcurl\solw}{\VSL(\Omega)},
\end{align*}
see, e.g., \cite[Sec. 3.5.3]{BookMonk}. 

\begin{remark}
If there holds $\solv\in\Hdivgamma{\Gamma}\cap\VSL_T(\Gamma)$ as well as $\solw\in\Hcurlgamma{\Gamma}\cap\VSL_T(\Gamma)$, the duality pairing $\dualpair{\solv}{\solw}{\Zdiv}{\Ycurl}$ coincides with $\SCP{\solv}{\solw}{\VSL_T(\Gamma)}$.
\end{remark}

According to Proposition~\ref{traceprop} there holds $\solu_T\in\Hcurlgamma{\Gamma}$ as long as we have $\solu\in\Hcurl$. For more regular vector fields $\solu$ we expect higher regularity of the tangential component $\solu_T$. The following lemma shows that this is indeed the case.

\begin{lemma}\label{tracelemma}
		Let $\ell\in\N_0$ and let $\Omega$ be a bounded $\CM{\ell+2}$-domain with boundary $\Gamma$. Moreover, let $\solu\in\VSHM{\ell+1}(\Omega)$ satisfy $\curl\solu\in\VSHM{\ell+1}(\Omega)$. Then, there holds $\solu_T\in\VSHM{\ell+1/2}_T(\Gamma)$ as well as $\curl_{\Gamma}\solu_T\in\SHM{\ell+1/2}_T(\Gamma)$ with 
	\begin{align*}
		\norm{\solu_T}{\VSHM{\ell+1/2}_T(\Gamma)}+\norm{\curl_{\Gamma}\solu_T}{\SHM{\ell+1/2}_T(\Gamma)}\leq C \left(\norm{\solu}{\VSHM{\ell+1}(\Omega)}+\norm{\curl\solu}{\VSHM{\ell+1}(\Omega)}\right)
	\end{align*}
	for a constant $C>0$ depending only on $\ell$ and $\Omega$.
\end{lemma}
\begin{fatproof}
		The fact that $\solu_T\in\VSHM{\ell+1/2}_T(\Gamma)$ follows from the trace theorem  and the fact that the outer unit normal field $\soln$ of a $\CM{\ell+2}$-domain is a $\VCM{\ell+1}$-vector field on $\Gamma$. Moreover, according to \eqref{surfacecurlconnect} there holds the equality $\curl_{\Gamma}\solu_T = \SCP{\curl\solu}{\soln}{\VSL(\Gamma)}$, thus the trace theorem and the assumed regularity of $\Gamma$ also prove $\curl_{\Gamma}\solu_T\in\SHM{\ell+1/2}_T(\Gamma)$ as well as the desired estimate.
\end{fatproof}

\subsection{Variational formulation of Maxwell's equations}\label{secweakform}

We still have to clarify the notion of a {\it weak solution} of Maxwell's equations \eqref{Maxwellorig}. To that end, let us assume that $\Gp=\geom$ is a $\CM{\ell}$-partition and that $\mu^{-1}$ and $\varepsilon$ are coercive $\VCMP{m}$-tensor fields for some $\ell\geq 2$ and $m\geq 1$, see Definition~\ref{partitiondef} and Definition~\ref{assumptioncoefs}. In addition, let the wavenumber $k\in\Co\setminus\{0\}$ and a right-hand side $\solf\in\VSL(\Omega)$ be given.

Moreover, depending on the imposed boundary condition, we make the following assumptions: 
\begin{itemize}
		\item In the case of impedance boundary conditions \eqref{Maxwellimpedance} we suppose that boundary data $\solgi\in\VSL_T(\Gamma)$ is given, and that $\zeta$ is a coercive $\VCM{m}$-tensor field acting in the tangent plane to $\Gamma$, see Definition~\ref{assumptionimpedance}. 
		\item In the case of natural boundary conditions \eqref{Maxwellnatural} we suppose that boundary data $\solgn\in\Hdivgamma{\Gamma}$ is given. 
\end{itemize}

Depending on the imposed boundary conditions, the energy space in which we look for a variational solution takes a different form. To account for that, we define
\begin{align*}
		&\HXI:=\{\solu\in\Hcurl\ |\ \solu_T\in\VSL(\Gamma)\}, \\
		&\HXE := \{\solu\in\Hcurl\ |\ \solu_T = 0\}.
\end{align*}

Depending on the boundary conditions, we call a vector field $\solu$ on $\Omega$ a weak solution of \eqref{Maxwellorig} if 
\begin{itemize}
		\item in the case of impedance boundary conditions \eqref{Maxwellimpedance}, there holds $\solu\in\HXI$ as well as
				\begin{align*}
					\SCP{\mu^{-1}\curl \solu}{\curl\solv}{\VSL(\Omega)}-k^2\SCP{\varepsilon\solu}{\solv}{\VSL(\Omega)}-ik\SCP{\zeta\solu_T}{\solv_T}{\VSL(\Gamma)}
					= \SCP{\solf}{\solv}{\VSL(\Omega)}+\SCP{\solgi}{\solv_T}{\VSL(\Gamma)}
				\end{align*}
				for all $\solv\in\HXI$, 

		\item in the case of natural boundary conditions \eqref{Maxwellnatural}, there holds $\solu\in\Hcurl$ as well as 
				\begin{align*}
					\SCP{\mu^{-1}\curl \solu}{\curl\solv}{\VSL(\Omega)}-k^2\SCP{\varepsilon\solu}{\solv}{\VSL(\Omega)} = 
					\SCP{\solf}{\solv}{\VSL(\Omega)}+\dualpair{\solgn}{\solv_T}{\Zdiv}{\Ycurl}				
			\end{align*}
				for all $\solv\in\Hcurl$,  
		\item in the case of essential boundary conditions \eqref{Maxwellessential}, there holds $\solu\in\HXE$ as well as 
			\begin{align*}
					\SCP{\mu^{-1}\curl \solu}{\curl\solv}{\VSL(\Omega)}-k^2\SCP{\varepsilon\solu}{\solv}{\VSL(\Omega)} = 
					\SCP{\solf}{\solv}{\VSL(\Omega)}				
			\end{align*}
			for all $\solv\in\HXE$.
\end{itemize}

%% file: regularity.tex
\section{Poisson transmission problems and Helmholtz decompositions}\label{transmissionproblems}

The aim of this section is to establish regularity shift properties of Poisson transmission problems with normal and tangential transmission conditions on the gradient of the solution $u$. Regularity of transmission problems has already been extensively discussed e.g. in \cite{TransmissionVogelius, TransmissionCaffarelli, TransmissionNirenberg, GrandLivre}. The canonical form of Poisson transmission problems is to have prescribed jumps of $u$ and its normal derivative. However, the problem that we consider only has prescribed jumps of the normal and tangential derivatives, not of $u$ itself. Therefore, an explicit proof of its regularity properties seems necessary.

After having discussed these transmission problems, we turn our attention to Helmholtz decompositions. Helmholtz decompositions of vector fields in $\VSL(\Omega)$ or $\Hcurl$ are a fundamental tool for the analysis of Maxwell's equations and have been thoroughly investigated in e.g. \cite{DecompositionAmrouche, BookCessenat, BookMonk}. We exploit recent results from \cite{CurlInverse, MaxwellImpedanceMelenk} to establish the existence of Helmholtz decompositions featuring piecewise regularity properties.


\subsection{Jump operators}\label{jumps}

In order to effectively discuss the effects of transmission conditions on the interface components $\interf_1,\ldots,\interf_r$ we define special jump operators.
Let $\ell\in\N_0$ and suppose that $\Gp=\geom$ is a $\CM{\ell+1}$-partition.
For a function $u\in\PSHM{\ell+1}(\Gp)$, the trace theorem guarantees that for $i=1,\ldots,n$ the trace $u_i\vert_{\partial \Gp_i}$ (where $u_i:=u\vert_{\Gp_i}$) is in $\SHM{\ell+1/2}(\partial \Gp_i)$. 

Moreover, if $\Gp$ is even a $\CM{\ell+2}$-partition and $\solu\in\PVSHM{\ell+1}(\Gp)$, the trace theorem asserts that for $i=1,\ldots,n$, the normal and tangential traces $\solu_i\vert_{\partial \Gp_i}\cdot\soln_{i}$ and $\solu_i\vert_{\partial \Gp_i}\times\soln_{i}$ (where $\solu_i := \solu\vert_{\Gp_i}$ and $\soln_i$ is the outer normal vector of the subdomain $\Gp_i$) are in $\SHM{\ell+1/2}(\partial \Gp_i)$ and $\VSHM{\ell+1/2}(\partial \Gp_i)$, respectively.

\medskip

Thus, for any $\CM{\ell+1}$-partition $\Gp$ we may define a jump operator $\jump{\ }\hskip -0.1cm:\PSHM{\ell+1}(\Gp)\rightarrow \SHM{\ell+1/2}(\interf)$ piecewise by 
\begin{align*}
		\jump{u} := u_{j_1}\vert_{\interf_j}-u_{j_2}\vert_{\interf_j},
\end{align*}
on every interface component $\interf_j = \partial\Gp_{j_1}\cap\partial\Gp_{j_2}$, where without loss of generality $j_1>j_2$, and where $u_{j_1} := u\vert_{\Gp_{j_1}}$ as well as $u_{j_2}:=u\vert_{\Gp_{j_2}}$.
Similarly, for any given $\CM{\ell+2}$-partition $\Gp$ we may define 
normal and tangential jump operators $\njump{\ }\hskip -0.1cm:\PVSHM{\ell+1}(\Gp)\rightarrow \SHM{\ell+1/2}(\interf)$ and 
$\tjump{\ }\hskip -0.1cm:\PVSHM{\ell+1}(\Gp)\rightarrow \VSHM{\ell+1/2}_T(\interf)$ piecewise by 
\begin{align*}
		\njump{\solu}\vert_{\interf_j} := \left(\solu_{j_1}\vert_{\interf_j}-\solu_{j_2}\vert_{\interf_j}\right)\cdot\soln_{j_1}\quad \rm{and}\quad \tjump{\solu}\vert_{\interf_j} := \soln_{j_1}\times\left[\left(\solu_{j_1}\vert_{\interf_j}-\solu_{j_2}\vert_{\interf_j}\right)\times\soln_{j_1}\right]
\end{align*} 
where $\soln_{j_1}$ is the outer normal vector to $\Gp_{j_1}$.

\bigskip 

For a $\CM{\ell+2}$-partition $\Gp$ and $z\in\SHM{\ell+3/2}(\interf)$ it is possible to construct a function $u\in\PSHM{\ell}(\Gp)$ satisfying $\norm{u}{\PSHM{\ell+2}(\Gp)}\leq C\norm{z}{\SHM{\ell+3/2}(\interf)}$ such that $u=0$ on $\Gamma$ and $\jump{u} = z$. We explain the procedure for a partition consisting of two subdomains $\Omega = \Gp_1\cup\Gp_2$, the general case follows the same ideas.
On a partition consisting of two subdomains, we
may define $u$ piecewise by setting $u\vert_{\Gp_1}:=0$ and $u\vert_{\Gp_2}:=v$, where $v$ is the solution of
\begin{align*}
\Delta v &= 0 \quad{\rm in}\ \Gp_2, \\
v &= z\quad{\rm on}\ \interf, \\
v &= 0\quad{\rm on}\ \partial \Gp_2\setminus\interf.
\end{align*}
By construction we have $\jump{u} = z$ as well as $\norm{u}{\PSHM{\ell+2}(\Gp)}\leq C\norm{z}{\SHM{\ell+3/2}(\interf)}$.

\medskip 

Similarly, for a $\CM{\ell+2}$-partition $\Gp$ and given $\solz\in\VSHM{\ell+1/2}_T(\interf)$ satisfying $\curl_{\interf}\solz = 0$, it is possible to construct a function $\solu\in\PVSHM{\ell+2}(\Gp)$ such that $\tjump{\solu} = \solz$. For the proof of this fact we need an auxiliary result.

\begin{lemma}\label{regularity:sgboundary}
		Let $\Gp$ be a $\CM{\ell+2}$-partition, and assume that the boundary $\Gamma$ consists of $n'\geq 1$ simply connected components. For a given coercive $\VCM{\ell+1}$-tensor field $\nu:\Omega\rightarrow\Co^{3\times 3}$, a right-hand side $f\in\SHM{\ell}(\Omega)$ and boundary data $\solg\in\VSHM{\ell+1/2}_T(\Gamma)$, consider the problem of finding $u\in\SHM{\ell+2}(\Omega)$ such that
	\begin{align*}
	-\diverg \nu\nabla u = &f\quad{\rm in}\ \Omega, \\
	\nabla_{\Gamma}u = & \solg\quad{\rm on}\ \Gamma.
	\end{align*}
	This problem has a solution $u\in\SHM{\ell+2}(\Omega)$ if and only if $\curl_{\Gamma}\left(\solg\right) = 0$, and the set of all solutions $u\in\SHM{\ell+2}(\Omega)$ forms an $n'$-dimensional affine space.
\end{lemma}

\begin{fatproof}
	Without loss of generality we assume that $\Gamma$ is simply connected; the general case follows from applying the subsequent arguments to every simply connected component of $\Gamma$.
	
	According to Proposition \ref{setting:Hodgedecomp} we have $\curl_{\Gamma}\nabla_{\Gamma}u = 0$, thus $\curl_{\Gamma}(\solg) = 0$ is a necessary condition. Suppose $\curl_{\Gamma}(\solg) = 0$. The idea is to use a similar trick as in \cite[Proof of Theorem 4.3]{MaxwellImpedanceMelenk}: According to Proposition \ref{setting:regularitysurfaceops} there holds  $\diverg_{\Gamma}(\solg)\in\SHM{\ell-1/2}(\Gamma)$, thus the zero mean-value solution $z$ of $\diverg_{\Gamma}\nabla_{\Gamma}z = \diverg_{\Gamma}(\solg)$ is an element of $\SHM{\ell+3/2}(\Gamma)$. We define $u$ as the solution of the Poisson problem
	\begin{align*}
			-\diverg\revision{\nu}\nabla u = &f\quad{\rm in}\ \Omega \\
	u = &z\quad{\rm on}\ \Gamma.
	\end{align*}
	Then we have $u\in\SHM{\ell+2}(\Omega)$ and \revision{from $\curl_{\Gamma}\left(\solg\right)=0$ and Proposition~\ref{setting:Hodgedecomp} we infer $\solg = \nabla_{\Gamma}\xi$ for some $\xi\in\SHM{1}(\Gamma)$ with zero mean. By construction of $u$ this leads to $\diverg_{\Gamma}\nabla_{\Gamma}u = \diverg_{\Gamma}\nabla_{\Gamma}z = \diverg_{\Gamma}\nabla_{\Gamma}\xi$ and by uniqueness of the zero-mean solution of this elliptic problem we infer $\nabla_{\Gamma}u = \nabla_{\Gamma}\xi = \solg$ on $\Gamma$.} 
	
	Finally, note that the difference $v := u_1-u_2$ of two solutions $u_1$ and $u_2$ satisfies $\diverg\nu\nabla v = 0$ in $\Omega$ and $\nabla_{\Gamma}v = 0$ on $\Gamma$. The latter implies that $v\vert_{\Gamma}$ is constant, thus, due to $\diverg\revision{\nu}\nabla v = 0$, the function $v$ is constant in $\Omega$. That is, the solution space forms a one-dimensional affine subspace of $\SHM{\ell+2}(\Omega)$.
\end{fatproof}

Let $\Gp$ be a $\CM{\ell+2}$-partition and suppose that $\solz\in\VSHM{\ell+1/2}_T(\interf)$ satisfies $\curl_{\Gamma}\solz = 0$. We claimed that it is possible to find $\solu\in\PVSHM{\ell+2}(\Gp)$ such that $\tjump{\solu} = \solz$. Again, we explain the procedure for a partition consisting of two subdomains $\Omega=\Gp_1\cup\Gp_2$. On such a partition 
we may define $\solu$  piecewise by $\solu\vert_{\Gp_1} := 0$ and $\solu\vert_{G_2}:=\nabla\varphi$, where $\varphi\in\SHM{\ell+2}(\Gp_2)$ is a solution of
\begin{align*}
\Delta \varphi &= 0\quad{\rm in}\ \Gp_2, \\
\nabla_{\interf}\varphi &= \solz\quad\rm{on}\ \interf, \\
\nabla_{\partial \Gp_2}\varphi &= 0\quad{\rm on}\ \partial \Gp_2\setminus\interf,
\end{align*}
which exists due to Lemma \ref{regularity:sgboundary}. It is straightforward to check that the piecewise defined function $\solu$ has the desired properties.

\input{regularity_1.tex}

\input{regularity_2.tex}

\section{Proof of Theorem~\ref{Mainresult3}}\label{finalsec}

We conclude this work by giving a proof of the third main result, Theorem~\ref{Mainresult3}. 

\bigskip

\begin{fatproofmod}{Theorem~\ref{Mainresult3}}
We split the proof into three steps. The first step proves the estimate corresponding to $\solu$ and the second step discusses the situation for $\curl\solu$. The third step then concludes the proof.	
	\smallskip 
	
	{\bf Step 1:} The first step is to show that for any $m\in\N_0$ with $m\leq \ell$ there holds $\solu\in\PVSHM{m+1}(\Gp)$ with the estimate
	\begin{align}\label{finiteregstep1}
			\begin{split}
				|k|^{-m-1}\norm{\solu}{\PVSHM{m+1}(\Gp)}\leq C_m&\left(|k|^{-m-1}\norm{\curl\solu}{\PVSHM{m}(\Gp)}+|k|^{-m}\norm{\solu}{\PVSHM{m}(\Gp)}+|k|^{-2}\left[\solF_{m,k}+\solG_{m,k}\right]\right)
			\end{split}
			\end{align}
	for a constant $C_m>0$ depending only on $m$, $\Gp$, $\mu^{-1}$, $\varepsilon$ and, if impedance boundary conditions are imposed, also on $\zeta$.
	
	We proceed by induction: For $m\in\N_0$ with $m\leq\ell$ assume that $\solu\in\PVSHM{m}(\Gp)$ with $\curl\solu\in\PVSHM{m}(\Gp)$. Note that for $m=0$ this assumption is certainly satisfied, since $\solu\in\Hcurl$. 
	
	Depending on the imposed boundary conditions, we distinguish between three cases.

	\medskip 

	{\bf Case 1:} Suppose that impedance boundary conditions \eqref{Maxwellimpedance} are imposed. By Theorem~\ref{Mainresult2} we have 
	\begin{align*}
		\norm{\solu}{\PVSHM{m+1}(\Gp)}\leq C\left(\norm{\curl\solu}{\PVSHM{m}(\Gp)}+\norm{\diverg\varepsilon\solu}{\PSHM{m}(\Omega)}+\norm{\solu_T}{\VSHM{m+1/2}_T(\Gamma)}\right),
	\end{align*}
	hence due to $-k^2\diverg\varepsilon\solu = \diverg\solf$ we get
	\begin{align}\label{impedancecase}
			\norm{\solu}{\PVSHM{m+1}(\Gp)}\leq C\left(\norm{\curl\solu}{\PVSHM{m}(\Gp)}+|k|^{-2}\norm{\diverg\solf}{\PSHM{m}(\Omega)}+\norm{\solu_T}{\VSHM{m+1/2}_T(\Gamma)}\right),
	\end{align}
	It remains to estimate $\norm{\solu_T}{\VSHM{m+1/2}(\Gamma)}$. According to Lemma~\ref{shiftresultsurface} we have 
	\begin{align*}
		\norm{\solu_T}{\VSHM{m+1/2}_T(\Gamma)}\leq C\left(\norm{\diverg_{\Gamma}\zeta\solu_T}{\SHM{m-1/2}(\Gamma)}+\norm{\curl_{\Gamma}\solu_T}{\SHM{m-1/2}(\Gamma)}\right), 
	\end{align*}
	and due to $\curl_{\Gamma}\solu_T = \curl\solu\cdot\soln$ on $\Gamma$, there holds $\norm{\curl_{\Gamma}\solu_T}{\revision{\SHM{m-1/2}(\Gamma)}}\leq C\norm{\curl\solu}{\PVSHM{m}(\Gp)}$. In addition, taking the surface divergence on the impedance boundary condition yields 
	\begin{align*}
		ik\diverg_{\Gamma}\zeta\solu_T&=\diverg_{\Gamma}(\mu^{-1}\curl\solu\times\soln)-\diverg_{\Gamma}\solgi = \curl\mu^{-1}\curl\solu\cdot\soln-\diverg_{\Gamma}\solgi \\
		& = \solf\cdot\soln+k^2\varepsilon\solu\cdot\soln-\diverg_{\Gamma}\solgi,
	\end{align*}
	and thus 
	\begin{align*}
			\norm{\diverg_{\Gamma}\zeta\solu_T}{\VSHM{m-1/2}(\Gamma)}\leq C\left(|k|\norm{\solu}{\PVSHM{m}(\Gp)}+\revision{|k|^{-1}\norm{\solf}{\PVSHM{m}(\Gp)}}+|k|^{-1}\norm{\diverg\solf}{\PSHM{m}(\Gp)}+|k|^{-1}|\norm{\diverg_{\Gamma}\solgi}{\SHM{m-1/2}(\Gamma)}\right). 
	\end{align*}
	In total we end up with 
	\begin{align*}
			\norm{\solu_T}{\VSHM{m+1/2}(\Gamma)}\leq C&\left(\norm{\curl\solu}{\PVSHM{m}(\Gp)}+|k|\norm{\solu}{\PVSHM{m}(\Gp)}+\revision{|k|^{-1}\norm{\solf}{\PVSHM{m}(\Gp)}}+|k|^{-1}\norm{\diverg\solf}{\PSHM{m}(\Gp)}\right.  \\ 
													  & \hskip 9cm \left. +|k|^{-1}\norm{\diverg_{\Gamma}\solgi}{\SHM{m-1/2}(\Gamma)}\right), 
	\end{align*}
	and together with \eqref{impedancecase} and a multiplication with $|k|^{-m-1}$ this yields \eqref{finiteregstep1} in the case of impedance boundary conditions.

	\medskip 

	{\bf Case 2:} Suppose that natural boundary conditions \eqref{Maxwellnatural} are imposed. By Theorem~\ref{Mainresult1} and the equality $-k^2\diverg\varepsilon\solu = \diverg\solf$ we have 
	\begin{align}\label{naturalcase1}
			\norm{\solu}{\PVSHM{m+1}(\Gp)}\leq C\left(\norm{\curl\solu}{\PVSHM{m}(\Gp)}+|k|^{-2}\norm{\diverg\solf}{\PSHM{m}(\Omega)}+\norm{\varepsilon\solu\cdot\soln}{\SHM{m+1/2}(\Gamma)}\right).
	\end{align}
	Moreover, 
	\begin{align*}
		\diverg_{\Gamman}\solgn = \diverg_{\Gamman}(\mu^{-1}\curl\solu\times\soln) = \curl\mu^{-1}\curl\solu\cdot\soln=\solf\cdot\soln+k^2\varepsilon\solu\cdot\soln,
	\end{align*}
	hence 
	\begin{align}\label{naturalcase2}
	\norm{\varepsilon\solu\cdot\soln}{\SHM{m+1/2}(\Gamma)}\leq |k|^{-2}\norm{\solf\cdot\soln\revision{-\diverg_{\Gamma}\solgn}}{\SHM{m+1/2}(\Gamma)}.
	\end{align}
	Combining \eqref{naturalcase1} and \eqref{naturalcase2} and multiplying with $|k|^{-m-1}$ yields \eqref{finiteregstep1}	in the case of natural boundary conditions.

	\medskip 

	{\bf Case 3:} Suppose that essential boundary conditions \eqref{Maxwellessential} are imposed. Due to Theorem~\ref{Mainresult2} as well as $-k^2\diverg\varepsilon\solu = \diverg\solf$ and $\solu_T=0$ on $\Gamma$ we have
	\begin{align*}
			\norm{\solu}{\PVSHM{m+1}(\Gp)}\leq C\left(\norm{\curl\solu}{\PVSHM{m}(\Gp)}+|k|^{-2}\norm{\diverg\solf}{\PSHM{m}(\Omega)}\right),
	\end{align*}
	from which a multiplication with $|k|^{-m-1}$ implies \eqref{finiteregstep1} in the case of essential boundary conditions. 


	\medskip 

	In total, this proves \eqref{finiteregstep1} for all three boundary conditions \eqref{Maxwellimpedance}, \eqref{Maxwellnatural} and \eqref{Maxwellessential} and thus concludes Step~1 of the proof.
	
	\bigskip 

	{\bf Step 2:} In this step we show that for all $m\in\N_0$ with $m\leq\ell$ we have $\curl\solu\in\PVSHM{m+1}(\Gp)$ with the estimate
	\begin{align}\label{finiteregstep2}
					|k|^{-m-2}\norm{\curl\solu}{\PVSHM{m+1}(\Gp)}\leq C_m\left(|k|^{-m-1}\norm{\curl\solu}{\PVSHM{m}(\Gp)}+|k|^{-m}\norm{\solu}{\PVSHM{m}(\Gp)}
																		 +|k|^{-2}\left[\solF_{m,k}+\solG_{m,k}\right]\right)
	\end{align}
	for a constant $C_m>0$ depending only on $m, \Gp, \varepsilon$, $\mu^{-1}$ and, if impedance boundary conditions are posed, also on $\zeta$.
	
	Again, we proceed by induction on $m$ and suppose $\solu\in\PVSHM{m}(\Gp)$, which for $m=0$ is certainly satisfied, and as before, we consider three cases.

	\medskip 

	{\bf Case 1:} Suppose that impedance boundary conditions \eqref{Maxwellimpedance} are prescribed. Due to the assumed coercivity and piecewise regularity of $\mu^{-1}$ we have 
\begin{align*}
\norm{\curl\solu}{\PVSHM{m+1}(\Gp)}\leq C\norm{\mu^{-1}\curl\solu}{\PVSHM{m+1}(\Gp)},
\end{align*}
thus Theorem~\ref{Mainresult2} and $\diverg\curl\solu = 0$ imply
	\begin{align*}
			\norm{\curl\solu}{\PVSHM{m+1}(\Gp)}&\leq C\left(\norm{\curl\mu^{-1}\curl\solu}{\PVSHM{m}(\Gp)}+\norm{(\mu^{-1}\curl\solu)_t}{\VSHM{m+1/2}(\Gamma)}\right).
	\end{align*}

    Maxwell's equations \eqref{Maxwellorig} and the assumed piecewise regularity of $\varepsilon$ yield
	\begin{align*}
		\norm{\curl\mu^{-1}\curl\solu}{\PVSHM{m}(\Gp)}\leq \norm{\solf}{\PVSHM{m}(\Gp)}+|k|^2\norm{\solu}{\PVSHM{m}(\Gp)}, 
	\end{align*}
	hence we obtain
	\begin{align}\label{impedancecasestep2}
			\norm{\curl\solu}{\PVSHM{m+1}(\Gp)}&\leq C\left(\norm{\solf}{\PVSHM{m}(\Gp)}+|k|^2\norm{\solu}{\PVSHM{m}(\Gp)}+\norm{(\mu^{-1}\curl\solu)_t}{\VSHM{m+1/2}(\Gamma)}\right).
	\end{align}
	Furthermore, due to Lemma~\ref{shiftresultsurface} there holds 
	\begin{align*}
			\norm{(\mu^{-1}\curl\solu)_t}{\VSHM{m+1/2}(\Gamma)}\leq C\left(\norm{\curl_{\Gamma}\zeta^{-1}(\mu^{-1}\curl\solu)_t}{\SHM{m-1/2}(\Gamma)}+\norm{\diverg_{\Gamma}(\mu^{-1}\curl\solu)_t}{\SHM{m-1/2}(\Gamma)}\right).
	\end{align*}
	We notice that 
	\begin{align*}
			\diverg_{\Gamma}(\mu^{-1}\curl\solu)_t = \curl\mu^{-1}\curl\solu\cdot\soln = \solf\cdot\soln+k^2\varepsilon\solu\cdot\soln,
	\end{align*}
	and the impedance condition \eqref{Maxwellimpedance} reads as $(\mu^{-1}\curl\solu)_t-ik\zeta\solu_T = \solg_T$, hence
	\begin{align*}
			\curl_{\Gamma}\zeta^{-1}(\mu^{-1}\curl\solu)_t = ik\curl_{\Gamma}\solu_T+\curl_{\Gamma}(\zeta^ {-1}\solgi) = ik\curl\solu\cdot\soln+\curl_{\Gamma}(\zeta^{-1}\solgi).
	\end{align*}
	All in all, we infer
	\begin{align*}
			\norm{(\mu^{-1}\curl\solu)_t}{\VSHM{m+1/2}(\Gamma)}\leq C\left(|k|\norm{\curl\solu}{\PVSHM{m}(\Gp)}+|k|^2\norm{\solu}{\PVSHM{m}(\Gp)}+|k|^m\left[\solF_{m,k}+\solG_{m,k}\right]\right).
	\end{align*}
	Together with \eqref{impedancecasestep2} and a multiplication with $|k|^{-m-2}$, this yields \eqref{finiteregstep2} in the case of impedance boundary conditions.
	
	\medskip

	{\bf Case 2:} Suppose that natural boundary conditions \eqref{Maxwellnatural} are prescribed. We notice that \eqref{impedancecasestep2} still holds true in this case, hence the boundary condition $(\mu^{-1}\curl\solu)_t = \solgn$ and a multiplication with $|k|^{-m-2}$ conclude the proof of \eqref{finiteregstep2} in the case of natural boundary conditions.
	
\medskip 

{\bf Case 3:} Suppose that essential boundary conditions \eqref{Maxwellessential} are prescribed. Analogously to the derivation of \eqref{impedancecasestep2} but with Theorem~\ref{Mainresult1} instead of Theorem~\ref{Mainresult2} we arrive at
	\begin{align*}
			\norm{\curl\solu}{\PVSHM{m+1}(\Gp)}&\leq C\left(\norm{\solf}{\PVSHM{m}(\Gp)}+|k|^2\norm{\solu}{\PVSHM{m}(\Gp)}+\norm{\curl\solu\cdot\soln}{\SHM{m+1/2}(\Gamma)}\right).
	\end{align*}
	From the boundary condition $\solu_T=0$ on $\Gamma$ we infer $\curl\solu\cdot\soln=0$ on $\Gamma$, hence a multiplication with $|k|^{-m-2}$ shows \eqref{finiteregstep2} in the case of essential boundary conditions.
	
	\bigskip 

	{\bf Step 3:} The third step is to combine \eqref{finiteregstep1} and \eqref{finiteregstep2} to conclude \eqref{finiteregularity}. For $m\in\N_0$ with $m\leq\ell$ we define the quantity
	\begin{align*}
		\vartheta_{m,k}:=|k|^{-m}\norm{\solu}{\PVSHM{m}(\Gp)}+|k|^{-m-1}\norm{\curl\solu}{\PVSHM{m}(\Gp)}.
	\end{align*}
	Combining \eqref{finiteregstep1} and \eqref{finiteregstep2} yields 
	\begin{align*}
			\vartheta_{m+1,k}\leq C_m\left(\vartheta_{m,k}+|k|^{-2}\left[\solF_{m,k}+\solG_{m,k}\right]\right),
	\end{align*}
	and by induction we infer 
	\begin{align*}
			\vartheta_{\revision{m}+1,k}\leq C_{\revision{m}}\left(\vartheta_{0,k}+|k|^{-2}\sum_{j=0}^{\revision{m}}\left[\solF_{j,k}+\solG_{j,k}\right]\right).
	\end{align*}
	By definition of $\vartheta_{0,k}$ we get 
	\begin{align*}
			\vartheta_{\revision{m}+1,k}\leq C_{\revision{m}}\left(\norm{\solu}{\VSL(\Omega)}+|k|^{-1}\norm{\curl\solu}{\VSL(\Omega)}+|k|^{-2}\sum_{j=0}^{\revision{m}}\left[\solF_{j,k}+\solG_{j,k}\right]\right)
	\end{align*}
	\revision{for all $m\in\{0,\ldots,\ell\}$. This is precisely \eqref{finiteregularity}, hence the proof of Theorem~\ref{Mainresult3} is complete.}
\end{fatproofmod}
\section*{Acknowledgements}
We gladly acknowledge financial support by Austrian Science Fund (FWF)
throughout the
special research program \textit{Taming complexity in PDE systems} (grant SFB F65,
\href{https://doi.org/10.55776/F65}{DOI:10.55776/F65}; JMM).

%% file: regularity_1.tex
\subsection{Piecewise Sobolev regularity of Poisson transmission problems}

The subsequent result establishes Sobolev regularity properties of solutions of Poisson transmission problems with prescribed jumps of the normal and tangential derivatives.

\begin{lemma}\label{regularity:sgtransmission}
		Let $\ell\in\N_0$, let $\Gp$ be a $\CM{\ell+2}$-partition and assume that $\nu$ is a coercive $\VCMP{\ell+1}$-tensor field in the sense of Definition~\ref{assumptioncoefs}.
		Suppose that $u\in\PSHM{1}(\Gp)$ satisfies
		\begin{align*}
		\begin{split}
				-\diverg\nu\nabla u = \diverg\solf&\hskip 0.5cm \rm{in}\ \Gp_1\cup\ldots\cup\Gp_n, \\
                \nu\nabla u\cdot\soln = f&\hskip 0.5cm \rm{on}\ \Gamma, \\
		\tjump{\nabla u} = \solg&\hskip 0.5cm \rm{on}\ \interf, \\
		\njump{\nu\nabla u} = h&\hskip 0.5cm \rm{on}\ \interf,
		\end{split}
		\end{align*}
		where $\solf\in\PVHdiv{\ell}$ and $f\in\SHM{\ell+1/2}(\Gamma)$, as well as $\solg\in\VSHM{\ell+1/2}_T(\interf)$ and $h\in\SHM{\ell+1/2}(\interf)$.
		
		Then, there holds $u\in\PSHM{\ell+2}(\Gp)$ and\hspace{0.1cm}\footnote{\revision{By a slight abuse of notation we write $\diverg\solf$ for the piecewise divergence of $\solf$.}} 
	\begin{align}\label{transmission}
	\norm{u}{\PSHM{\ell+2}(\Gp)}\leq C\left(\norm{\diverg\solf}{\PSHM{\ell}(\Gp)}+\norm{f}{\SHM{\ell+1/2}(\Gamma)}+\norm{\solg}{\VSHM{\ell+1/2}_T(\interf)}+\norm{h}{\SHM{\ell+1/2}(\interf)}+\norm{u}{\SL(\Omega)}\right),
	\end{align}
	where the constant $C>0$ depends only on $\Gp$, $\nu$ and $\ell$. 
	
	Suppose that in addition we have $\SCP{u}{1}{\SL(\Omega)}=0$ and $\SCP{\jump{u}}{1}{\SL(\interf)} = 0$. Then, \eqref{transmission} can be improved to 
	\begin{align}\label{transmission2}
	\norm{u}{\PSHM{\ell+2}(\Gp)}\leq C\left(\norm{\diverg\solf}{\PSHM{\ell}(\Gp)}+\norm{f}{\SHM{\ell+1/2}(\Gamma)}+\norm{\solg}{\VSHM{\ell+1/2}_T(\interf)}+\norm{h}{\SHM{\ell+1/2}(\interf)}\right),
	\end{align}
	where the constant $C>0$ again depends only on $\Gp$, $\nu$ and $\ell$.
\end{lemma}

\begin{fatproof}
	The proof is divided into three steps. The first step deals with \eqref{transmission2}, and the second and third step prove \eqref{transmission}.
	
	{\bf Step 1:}
	We first prove \eqref{transmission2} provided that $\SCP{u}{1}{\SL(\Omega)}=0$, and $\SCP{\jump{u}}{1}{\SL(\interf)} = 0$. We start connecting the tangential jump $\solg:=\tjump{\nabla u}$ to $\jump{u}$.
	To that end, we note that  there holds $\tjump{\nabla u} = \nabla_{\interf}\jump{u}$, and therefore $\diverg_{\interf}\tjump{\nabla u} = \diverg_{\interf}\nabla_{\interf}\jump{u}\in\SHM{\ell-1/2}(\interf)$.
	Hence, due to elliptic regularity and $\SCP{\jump{u}}{1}{\SL(\interf)} = 0$ we have $\jump{u}\in\SHM{\ell+3/2}(\interf)$ with $\norm{\jump{u}}{\SHM{\ell+3/2}(\interf)}\leq C\norm{\solg}{\VSHM{\ell+1/2}_T(\interf)}$.
	
	\medskip 
	
	As discussed in Section \ref{jumps}, there exists a function $w\in\PSHM{\ell+2}(\Gp)$ with $\jump{w} = \jump{u}$ and which satisfies the inequality
	\begin{align}\label{transmissiontmp0}
	\norm{w}{\PSHM{\ell+2}(\Gp)}\leq C\norm{\jump{u}}{\SHM{\ell+3/2}(\interf)}\leq C\norm{\solg}{\VSHM{\ell+1/2}_T(\interf)}.
	\end{align}
	We may redefine $w$ to $w-\frac{1}{|\Omega|}\SCP{w}{1}{\SL(\Omega)}$ without changing the property $\jump{w} = \jump{u}$ and we note that $\SCP{w}{1}{\SL(\Omega)}\leq C\norm{w}{\SL(\Omega)}$, thus the redefined $w:=w-\frac{1}{|\Omega|}\SCP{w}{1}{\SL(\Omega)}$ still satisfies \eqref{transmissiontmp0}. That is, we may without loss of generality assume $\SCP{w}{1}{\SL(\Omega)} = 0$.
	
	Let $\diverg\widetilde{\solf} := \diverg\solf+\diverg\nu\nabla w\in\PSHM{\ell}(\Gp)$ and $\widetilde{h} := \njump{\nu\nabla u-\nu\nabla w}\in\SHM{\ell+1/2}(\interf)$. Then, $v := u-w$ satisfies $-\diverg\nu\nabla v = \diverg\widetilde{\solf}$ on \revision{all subdomains $\Gp_j$}, as well as $\njump{\nu\nabla v} = \widetilde{h}$ and $\jump{v} = 0$, and $\nu\nabla v\cdot \soln = f-\nu\nabla w\cdot \soln$ on the boundary $\Gamma$. Moreover, we note that $\SCP{v}{1}{\SL(\Omega)} = 0$.
	
	\medskip 
	
	The equation $\jump{v} = 0$ implies $v\in\SHM{1}(\Omega)$, and partial integration on \revision{all subdomains $\Gp_j$} yields that $v$ satisfies 
	\begin{align}\label{transmissiontmp}
	\SCP{\nu\nabla v}{\nabla\xi}{\VSL(\Omega)} = \SCP{\widetilde{h}}{\xi}{\SL(\interf)}+\SCP{\widetilde{f}}{\xi}{\SL(\Gamma)}+\sum_{i=1}^n\SCP{\diverg\widetilde{\solf}}{\xi}{\SL(\Gp_i)}
	\end{align}
	for all $\xi\in\SHM{1}(\Omega)$, where $\widetilde{f} := f-\nu\nabla w\cdot \soln$. According to \cite[Theorem 5.3.8]{GrandLivre}, the function $v$ belongs to the space $\PSHM{\ell+2}(\Gp)$ and 
	\begin{align}\label{transmissiontmp2}
	\norm{v}{\PSHM{\ell+2}(\Gp)}\leq C\left(\norm{\diverg\widetilde{\solf}}{\PSHM{\ell}(\Gp)}+\norm{\widetilde{h}}{\SHM{\ell+1/2}(\interf)}+\norm{\widetilde{f}}{\SHM{\ell+1/2}(\Gamma)}+\norm{v}{\SHM{1}(\Omega)}\right).
	\end{align}

	\revision{By choosing $\xi=v$ in \eqref{transmissiontmp}, a Poincar\'{e} inequality (which is available due to $\SCP{v}{1}{\SL(\Omega)}=0$) and some straightforward estimates we obtain
			\begin{align*}
					\norm{\nabla v}{\VSL(\Omega)}\leq C\norm{\widetilde{h}}{\SL(\interf)}+\norm{\widetilde{f}}{\SL(\Gamma)}+\norm{\diverg\widetilde{\solf}}{\SL(\Omega)}
			\end{align*}
			for a constant $C>0$ depending only on $\nu$ and $\Omega$, where we slightly abused notation and wrote $\diverg\widetilde{\solf}$ for the piecewise divergence of $\widetilde{\solf}$.Exploiting the Poincar\'{e} inequality a second time yields 
\begin{align*}
		\norm{v}{\SHM{1}(\Omega)}\leq C\norm{\widetilde{h}}{\SL(\interf)}+\norm{\widetilde{f}}{\SL(\Gamma)}+\norm{\diverg\widetilde{\solf}}{\SL(\Omega)},
\end{align*}}
	hence \eqref{transmissiontmp2} can be improved to 
	\begin{align}\label{transmissiontmp3}
	\norm{v}{\PSHM{\ell+2}(\Gp)}\leq C\left(\norm{\diverg\widetilde{\solf}}{\PSHM{\ell}(\Gp)}+\norm{\widetilde{h}}{\SHM{\ell+1/2}(\interf)}+\norm{\widetilde{f}}{\SHM{\ell+1/2}(\Gamma)}\right).
	\end{align}
	Together with $v:=u-w$, the definitions of $\diverg\widetilde{\solf}$, $\widetilde{f}$ and $\widetilde{h}$, as well as \eqref{transmissiontmp0}, the inequality \eqref{transmissiontmp3} implies \eqref{transmission2}.
	
	\medskip 
	
	{\bf Step 2:} We assume that neither $\SCP{u}{1}{\SL(\Omega)} = 0$ nor $\SCP{\jump{u}}{1}{\SL(\interf)} = 0$ is satisfied anymore. In this step we prove the inequality 
	\begin{align}\label{step2}
	\norm{u}{\PSHM{\ell+2}(\Gp)}\leq C\left(\norm{\diverg\solf}{\PSHM{\ell}(\Gp)}+\norm{f}{\SHM{\ell+1/2}(\Gamma)}+\norm{\solg}{\VSHM{\ell+1/2}_T(\interf)}+\norm{h}{\SHM{\ell+1/2}(\interf)}+\norm{u}{\PSHM{1}(\Gp)}\right).
	\end{align}
	The proof of \eqref{step2} follows the same lines as the proof of \eqref{transmission2}, but instead of \eqref{transmissiontmp0}, the lifting $w\in\PSHM{\ell+2}(\Gp)$ of $\jump{u}$ satisfies 
	\begin{align*}
	\norm{w}{\PSHM{\ell+2}(\Gp)}\leq C\norm{\jump{u}}{\SHM{\ell+3/2}(\interf)}\leq C\left(\norm{\solg}{\VSHM{\ell+1/2}_T(\interf)}+\norm{\jump{u}}{\SL(\interf)}\right) \leq C\left(\norm{\solg}{\VSHM{\ell+1/2}_T(\interf)}+\norm{u}{\PSHM{1}(\Gp)}\right).
	\end{align*}
	Analogously as in step 1, one proves \eqref{transmissiontmp2} for $v:=u-w$. Then, \eqref{step2} follows from \eqref{transmissiontmp2}, the inequality $\norm{v}{\SHM{1}(\Omega)}\leq \norm{u}{\PSHM{1}(\Gp)}+\norm{w}{\PSHM{1}(\Omega)}$ and $\norm{w}{\PSHM{\ell+2}(\Gp)} \leq C\left(\norm{\solg}{\VSHM{\ell+1/2}_T(\interf)}+\norm{u}{\PSHM{1}(\Omega)}\right)$. 
	
	\medskip 
	
	{\bf Step 3:} The third and final step of the proof is to conclude \eqref{transmission}. From \eqref{step2} we see that it suffices to bound $\norm{\nabla u}{\VSL(\Omega)}$ in terms of $\diverg\solf$, $f$, $\solg$ and $h$ in suitable norms. By adding a suitable piecewise constant function to $u$ we can construct a function  $\widetilde{u}\in\PSHM{\ell+2}(\Gp)$ with $\nabla\widetilde{u} = \nabla u$ on every subdomain $\Gp_i$ and $\SCP{\widetilde{u}}{1}{\SL(\Omega)} = 0$, as well as $\SCP{\revision{\jump{\widetilde{u}}}}{1}{\SL(\interf)} = 0$. 
	According to step 1, the function $\widetilde{u}$ satisfies \eqref{transmission2}, which implies 
	\begin{align*}
	\norm{\nabla u}{\VSL(\Omega)}\leq C\left(\sum_{i=1}^n\norm{\diverg\solf}{\SL(\Gp_i)}+\norm{f}{\SHM{1/2}(\Gamma)}+\norm{\solg}{\VSHM{1/2}_T(\interf)}+\norm{h}{\SHM{1/2}(\interf)}\right).
	\end{align*}
	Together with \eqref{step2} this leads to \eqref{transmission}, which concludes the proof.
	
\end{fatproof}

In Lemma \ref{regularity:sgtransmission} we only considered Neumann boundary conditions on $\Gamma$. We note that it is also possible to consider Dirichlet boundary conditions, as the next lemma shows. 

\begin{lemma}\label{regularity:sgtransmissiondir}
		Let $\ell\in\N_0$, let $\Gp$ be a $\CM{\ell+2}$-partition and assume that $\nu$ is a coercive $\VCMP{\ell+1}$-tensor field in the sense of Definition~\ref{assumptioncoefs}.
		Suppose that $u\in\PSHM{1}(\Gp)$ satisfies
		\begin{align*}
		\begin{split}
				-\diverg\nu\nabla u = \diverg\solf&\hskip 0.5cm \rm{in}\ \Gp_1\cup\ldots\cup\Gp_n, \\
                u = f&\hskip 0.5cm \rm{on}\ \Gamma, \\
		\tjump{\nabla u} = \solg&\hskip 0.5cm \rm{on}\ \interf, \\
		\njump{\nu\nabla u} = h&\hskip 0.5cm \rm{on}\ \interf,
		\end{split}
		\end{align*}
		where $\solf\in\PVHdiv{\ell}$ and $f\in\SHM{\ell+3/2}(\Gamma)$, as well as $\solg\in\VSHM{\ell+1/2}_T(\interf)$ and $h\in\SHM{\ell+1/2}(\interf)$.
		
	Then, there holds $u\in\PSHM{\ell+2}(\Gp)$ and 
	\begin{align}\label{transmissiondir}
	\norm{u}{\PSHM{\ell+2}(\Gp)}\leq C\left(\norm{\diverg\solf}{\PSHM{\ell}(\Gp)}+\norm{f}{\SHM{\ell+3/2}(\Gamma)}+\norm{\solg}{\VSHM{\ell+1/2}_T(\interf)}+\norm{h}{\SHM{\ell+1/2}(\interf)}+\norm{u}{\PSHM{1}(\Gp)}\right),
	\end{align}
	where the constant $C>0$ depends only on $\Gp$, $\nu$ and $\ell$.

	If, in addition, $\SCP{\jump{u}}{1}{\SL(\interf)} = 0$, then \eqref{transmissiondir} can be improved to 
	\begin{align*}
	\norm{u}{\PSHM{\ell+2}(\Gp)}\leq C\left(\norm{\diverg\solf}{\PSHM{\ell}(\Gp)}+\norm{f}{\SHM{\ell+3/2}(\Gamma)}+\norm{\solg}{\VSHM{\ell+1/2}_T(\interf)}+\norm{h}{\SHM{\ell+1/2}(\interf)}\right),
	\end{align*}
	where the constant $C>0$ again depends only on $\Gp$, $\nu$ and $\ell$.
\end{lemma}
Since the proof of Lemma \ref{regularity:sgtransmissiondir} differs from the proof of Lemma \ref{regularity:sgtransmission} only in minor details, we omit it.

\begin{remark}
	Note that the right-hand side in \eqref{transmissiondir} involves the piecewise $\SHM{1}$-norm of $u$. This is in contrast to the corresponding inequality \eqref{transmission} in the case of Neumann boundary conditions where the right-hand side involves only the $\SL$-norm of $u$. Proving the inequality \eqref{transmissiondir} with $\norm{u}{\SL(\Omega)}$ instead of $\norm{u}{\PSHM{1}(\Gp)}$ might be possible, but seems to be more involved and is unnecessary for this work.
\end{remark}

\subsection{Helmholtz decompositions with piecewise regularity properties}
We turn our attention to the piecewise regularity properties of Helmholtz decompositions.
The subsequent two results state the existence of Helmholtz decompositions with piecewise regularity and orthogonality properties.

\begin{lemma}\label{reg:helmholtze}
		Let $\Gp$ be $\CM{2}$-partition and suppose that $\nu$ is a coercive $\VCMP{1}$-tensor field in the sense of Definition~\ref{assumptioncoefs}.
	 Furthermore, let $\solu\in\Hcurl$ be given. 
	\begin{itemize}
		\item[(i)] If $\SCP{\nu\solu}{\nabla\xi}{\VSL(\Omega)} = 0$ for all $\xi\in\SHM{1}(\Omega)$, then $\solu\in\PVSHM{1}(\Gp)$ and 
		\begin{align}\label{Helmholtzeineq1}
		\norm{\solu}{\PVSHM{1}(\Gp)}\leq C\norm{\curl\solu}{\VSL(\Omega)}.
		\end{align}
		\item[(ii)] There exists a decomposition $\solu = \solv+\nabla\varphi$ with 
		\begin{align}\label{Helmholtzeineq2}
		\norm{\solv}{\PVSHM{1}(\Gp)}\leq C\norm{\curl\solu}{\VSL(\Omega)},\quad \norm{\varphi}{\SHM{1}(\Omega)}\leq C\norm{\solu}{\VSL(\Omega)},
		\end{align}
		and $\SCP{\nu\solv}{\nabla\xi}{\VSL(\Omega)} = 0$ for all $\xi\in\SHM{1}(\Omega)$.
	\end{itemize}
	In both statements, the constant $C>0$ depends only on $\Gp$ and $\nu$. Moreover, if  $\nu\in\VCM{1}(\overline{\Omega})$, then \eqref{Helmholtzeineq1} and \eqref{Helmholtzeineq2} hold for $\norm{\solu}{\VSHM{1}(\Omega)}$ and $\norm{\solv}{\VSHM{1}(\Omega)}$ instead of $\norm{\solu}{\PVSHM{1}(\Gp)}$ and $\norm{\solv}{\PVSHM{1}(\Gp)}$, respectively.
\end{lemma}

\begin{fatproof}
		According to \cite[Lemma 2.7]{MaxwellImpedanceMelenk}, there exists a decomposition $\solu = \solv+\nabla\revision{\varphi}$, where $\solv$ and $\revision{\varphi}$ satisfy $\norm{\solv}{\VSHM{1}(\Omega)}\leq C\norm{\curl\solu}{\VSL(\Omega)}$ and $\norm{\revision{\varphi}}{\SHM{1}(\Omega)}\leq C\norm{\solu}{\Hcurl}$. 
		Due to $\nu\solu$ being orthogonal to all gradient fields, there holds $-\diverg\nu\nabla \revision{\varphi} = \diverg\nu\solv$ on every subdomain $\Gp_i$, as well as $\nu\nabla\revision{\varphi}\cdot\soln = -\nu\solv\cdot\soln$ on the boundary $\Gamma$. 
		Moreover, due to $\solv\in\VSHM{1}(\Omega)$ and $\nu\solu\in\Hdiv$ we have $\njump{\nabla\revision{\varphi}}=-\njump{\nu\solv}$ as well as $\tjump{\nabla\revision{\varphi}}=0$. Without loss of generality we may assume $\SCP{\revision{\varphi}}{1}{\SL(\Omega)} = 0$, and according to Lemma \ref{regularity:sgtransmission} the function $\revision{\varphi}$ satisfies 
	\begin{align*}
			\norm{\revision{\varphi}}{\PSHM{2}(\Gp)}\leq C\norm{\solv}{\SHM{1}(\Omega)}\leq C\norm{\curl\solu}{\VSL(\Omega)},
	\end{align*}
	where we used that $\revision{\varphi}\in\SHM{1}(\Omega)$ implies $\SCP{\jump{\revision{\varphi}}}{1}{\SL(\interf)}=0$. This proves $(i)$.
	
	In order to prove $(ii)$,
	let $\solu\in\Hcurl$ be given and let $\revision{\varphi}\in\SHM{1}(\Omega)$ satisfy $\SCP{\revision{\varphi}}{1}{\SL{\Omega}}=0$ and
	\begin{align*}
	\SCP{\nu\nabla\varphi}{\nabla\xi}{\VSL(\Omega)}= \SCP{\nu\solu}{\nabla\xi}{\VSL(\Omega)} 
\end{align*}
for all $\xi\in\SHM{1}(\Omega)$. Note that the Lax-Milgram lemma asserts the existence of $\revision{\varphi}$. 
	By construction, we have $\norm{\nabla\varphi}{\VSL(\Omega)}\leq C\norm{\solu}{\VSL(\Omega)}$ and a Poincar\'{e} inequality shows $\norm{\varphi}{\SHM{1}(\Omega)}\leq C \norm{\solu}{\VSL(\Omega)}$. Furthermore, we observe that $\solv :=\solu-\nabla\varphi\in\Hcurl$ satisfies $\SCP{\nu\solv}{\nabla\xi}{\VSL(\Omega)}=0$ for all $\xi\in\VSHM{1}(\Omega)$, and according to statement $(i)$ there holds $\norm{\solv}{\PVSHM{1}(\Gp)}\leq C\norm{\curl\solv}{\VSL(\Omega)} = C\norm{\curl\solu}{\VSL(\Omega)}$. This proves $(ii)$.
	
	If $\nu$ is continuously differentiable on $\overline{\Omega}$, then the proof follows the lines, except for the fact that the equation $-\diverg\nu\nabla\revision{\varphi} = \diverg\nu\solv$ in $\Omega$ with $\nu\nabla\revision{\varphi}\cdot\soln = -\nu\solv\cdot\soln$ on $\Gamma$ is a standard Neumann problem for $\revision{\varphi}$.
\end{fatproof}

The following result can be seen as a dual formulation of Lemma \ref{reg:helmholtze}.

\begin{lemma}\label{reg:helmholtzdiv}
		Let $\Gp$ be a $\CM{2}$-partition and assume that $\nu$ is a coercive $\VCMP{1}$-tensor field in the sense of Definition~\ref{assumptioncoefs}. 
      Furthermore, let $\solu$ satisfy $\nu\solu\in\Hdiv$.
	\begin{itemize}
		\item[(i)] If $\SCP{\solu}{\curl\solw}{\VSL(\Omega)}=0$ for all $\solw\in\Hcurl$, then $\solu\in\PVSHM{1}(\Gp)$ and
		\begin{align}\label{helmholtzdiveq1}
			\norm{\solu}{\PVSHM{1}(\Gp)}\leq C\norm{\diverg\nu\solu}{\SL(\Omega)}.
		\end{align}
		\item[(ii)] There exists a decomposition $\nu\solu = \nu\solv+\curl\solz$ with
		\begin{align}\label{helmholtzdiveq2}
			\norm{\solv}{\PVSHM{1}}\leq C\norm{\diverg\nu\solu}{\SL(\Omega)},\quad \norm{\solz}{\VSHM{1}(\Omega)}\leq \norm{\solu}{\VSL(\Omega)},
		\end{align}
		and $\SCP{\solv}{\curl\solw}{\VSL(\Omega)} = 0$ for all $\solw\in\Hcurl$.
	\end{itemize}
	In both statements, the constant $C>0$ depends only on $\Gp$ and $\nu$. Moreover, if $\nu\in\VCM{1}(\overline{\Omega})$, then \eqref{helmholtzdiveq1} and \eqref{helmholtzdiveq2} hold for $\norm{\solu}{\VSHM{1}(\Omega)}$ and $\norm{\solv}{\VSHM{1}(\Omega)}$, instead of $\norm{\solu}{\PVSHM{1}(\Gp)}$ and $\norm{\solv}{\PVSHM{1}(\Gp)}$, respectively.
\end{lemma}

\begin{fatproof}
	According to \cite[Lemma 3.27, (3.60)]{BookMonk}, the equation $\SCP{\solu}{\curl\solw}{\VSL(\Omega)}=0$ for all $\solw\in\Hcurl$ implies that $\solu = \nabla\xi$ for a function $\xi\in\SHM{1}_0(\Omega)$. Moreover, $\diverg\nu\nabla\xi = \diverg\nu\solu$, and therefore, according to Lemma \ref{regularity:sgtransmissiondir}, $\xi\in\PSHM{2}(\Gp)\cap\SHM{1}_0(\Omega)$ with $\norm{\xi}{\PSHM{2}(\Gp)}\leq C\norm{\diverg\nu\solu}{\SL(\Omega)}$.
	If $\nu$ is continuously differentiable on $\overline{\Omega}$, then $\diverg\nu\nabla\xi = \diverg\nu\solu$ implies $\norm{\xi}{\SHM{2}(\Omega)}\leq C\norm{\diverg\nu\solu}{\SL(\Omega)}$. This proves (i).
	
	In order to prove (ii), let $\solz$ satisfy $\SCP{\solz}{\nabla\xi}{\VSL{\Omega}}=0$ for all $\xi\in\SHM{1}(\Omega)$ and 
	\begin{align*}
			\revision{a(\solz,\solw):=}\SCP{\nu^{-1}\curl\solz}{\curl\solw}{\VSL(\Omega)} = \SCP{\solu}{\curl\solw}{\VSL(\Omega)}
	\end{align*}
	for all $\solw\in\Hcurl$. \revision{Note that according to Lemma~\ref{reg:helmholtze}, the sesquilinear for $a(\cdot,\cdot)$ is coercive on the space of vector fields $\solv\in\Hcurl$ that are orthogonal to all gradient fields, hence the Lax-Milgram lemma asserts the existence of $\solz$.} Moreover, we have $\norm{\curl\solz}{\VSL(\Omega)}\leq C\norm{\solu}{\VSL(\Omega)}$ and exploiting Lemma \ref{reg:helmholtze} leads to
	\begin{align*}
		\norm{\solz}{\VSHM{1}(\Omega)}\leq\norm{\curl\solz}{\VSL(\Omega)}\leq C\norm{\solu}{\VSL(\Omega)}.
	\end{align*}
	
	Moreover, $\solv := \solu-\nu^{-1}\curl\solz$ satisfies $\SCP{\solv}{\curl\solw}{\VSL(\Omega)} = 0$ for all $\solw\in\Hcurl$ and $\diverg\nu\solv = \diverg\nu\solu$, thus according to (i) there holds $\norm{\solv}{\PVSHM{1}\revision{(\Gp)}}\leq C\norm{\diverg\nu\solu}{\SL(\Omega)}$. 
	If $\nu$ is continuously differentiable on $\overline{\Omega}$, then (ii) is proved analogously, the only change is that the equations $\SCP{\solv}{\curl\solw}{\VSL(\Omega)} = 0$ for all $\solw\in\Hcurl$ and $\diverg\nu\solv = \diverg\nu\solu$ then imply $\norm{\solv}{\VSHM{1}\revision{(\Omega)}}\leq C\norm{\diverg\nu\solu}{\SL(\Omega)}$.
\end{fatproof}

%% file: regularity_2.tex
\section{Proof of Theorem~\ref{Mainresult1} and Theorem~\ref{Mainresult2}}\label{sec:regdecomp}

In this section we provide proofs for the first two main results of this work, namely Theorem~\ref{Mainresult1} and Theorem~\ref{Mainresult2}. 
The following proposition \cite[Lemma~2.6]{MaxwellImpedanceMelenk} has its origins in the seminal work \cite{CurlInverse} and is fundamental for our purposes. 

\begin{proposition}\label{regularity:costabel}
	Let $\Omega\subseteq \R^3$ be a bounded Lipschitz domain. There exist pseudodifferential operators $\Ro{\Omega}$, $\Rt{\Omega}$ of order $-1$ and $\Ko{\Omega}$, $\Kt{\Omega}$ of order $-\infty$ on $\R^3$ which for $m\in\Z$ have the mapping properties 
\begin{align*}
		\Ro{\Omega}:\VSHM{-m}(\Omega)&\rightarrow \SHM{1-m}(\Omega), \\
		\Rt{\Omega}:\VSHM{-m}(\Omega)&\rightarrow \VSHM{1-m}(\Omega), \\
		\Ko{\Omega},\Kt{\Omega}:\VSHM{m}(\Omega)&\rightarrow \left(C^{\infty}(\overline{\Omega})\right)^3.
\end{align*}
	 Moreover, for any $\solu\in\VSHM{m}(\Omega)\cap\Hcurlm{m}$ there holds
	\begin{align*}
	\solu = \nabla \Ro{\Omega}(\solu-\Rt{\Omega}(\curl \solu))+\Rt{\Omega}(\curl \solu)+\Ko{\Omega}\solu.
	\end{align*}
	In addition, $\Rt{\Omega}$ and $\Kt{\Omega}$ satisfy
	\begin{align*}
	\curl \Rt{\Omega} \solu = \solu-\Kt{\Omega}\solu
	\end{align*}
	for all $\solu\in\VSHM{m}(\Omega)\cap\Hcurlm{m}$ with $\diverg\solu=0$ on $\Omega$.
\end{proposition}

In essence, $\Rt{\Omega}$ is a right-inverse to the curl-operator; it plays a fundamental role in the subsequent proof of Theorem~\ref{Mainresult1}.

\medskip 

\begin{fatproofmod}{Theorem~\ref{Mainresult1}}
	 According to Lemma \ref{reg:helmholtze}, there exists a decomposition $\solu = \solv+\nabla\phi$ with  
	\begin{align*}
			\norm{\solv}{\PVSHM{1}(\Gp)}\leq C\norm{\curl\solu}{\VSL(\Omega)},\quad \norm{\revision{\phi}}{\SHM{1}(\Omega)}\leq C\norm{\solu}{\VSL(\Omega)},
	\end{align*}
	and $\SCP{\nu\solv}{\nabla\xi}{\VSL(\Omega)} = 0$ for all $\xi\in\SHM{1}(\Omega)$. 
	Together with $\nu\solu\cdot\soln = h$ on $\Gamma$ this implies $\diverg\nu\nabla\phi = \diverg\nu\solu$ on every subdomain $\Gp_i$  and $\nu\nabla\phi\cdot \soln = h$ on $\Gamma$. Moreover, we have $\njump{\nu\nabla\phi} = 0$, and there holds $\tjump{\nabla\phi} = 0$ on $\interf$ as well as $\SCP{\jump{\phi}}{1}{\SL(\interf)}=0$. Without loss of generality we may assume $\SCP{\phi}{1}{\SL(\Omega)}=0$, hence Lemma \ref{regularity:sgtransmission} yields
	\begin{align*}
			\norm{\phi}{\SHM{1}(\Omega)}+\norm{\nabla\phi}{\PVSHM{\ell+1}(\Gp)}\leq\norm{\phi}{\PSHM{\ell+2}}\leq C\left(\norm{\diverg\nu\solu}{\PVSHM{\ell}(\Gp)}+\norm{h}{\SHM{\ell+1/2}(\Gamma)}\right).
	\end{align*}
	In order to finish the proof, it suffices to show that $\solv$ satisfies 
	\begin{align}\label{sregtmp}
	\norm{\solv}{\PVSHM{\ell+1}(\Gp)}\leq C\norm{\curl \solu}{\PVSHM{\ell}(\Gp)},
	\end{align}
	which is done below by induction with respect to $\ell$.
	
	\medskip 
	
	We notice that for $\ell = 0$, the estimate \eqref{sregtmp} follows from Lemma \ref{reg:helmholtze}.
	 Assume that \eqref{sregtmp} holds for some $\ell\in\N_0$. We show that it is true for $\ell+1$ as well. For $i=1,\ldots,n$ we define $\solv_i:=\solv\vert_{\Gp_i}$ and apply Proposition~\ref{regularity:costabel} on every subdomain $\Gp_i$ to obtain
	\begin{align}\label{regularity:tmp2}
		\solv_i = \solz_i+\nabla\psi_i,
	\end{align}
	where $\solz_i:=\Rt{\Gp_i}(\curl\solv_i)+\Ko{\Gp_i}\solv_i$ and $\nabla\psi_i:= \nabla\Ro{\Gp_i}(\solv_i-\Rt{\Gp_i}(\curl \solv_i))$. The mapping properties of the operators $\Rt{\Gp_i}$ and $\Ko{\Gp_i}$, the equation $\curl\solu = \curl\solv$ and \eqref{sregtmp} imply 
	\begin{align}\label{reg:zinequality}
		\norm{\solz_i}{\PVSHM{\ell+2}(\Gp)}\leq C\left(\norm{\solv}{\PVSHM{\ell+1}(\Gp)}+\norm{\curl\solv}{\PVSHM{\ell+1}(\Gp)}\right)
		\leq C\left(\norm{\curl \solu}{\PVSHM{\ell+1}(\Gp)}\right).
	\end{align}
	
	We define $\solz\in\PVSHM{\ell+2}(\Gp)$ and $\psi\in\PSHM{\ell+2}(\Gp)$ piecewise by $\solz\vert_{\Gp_i}:=\solz_i$ and $\psi\vert_{\Gp_i}:=\phi_i$. Then, the decomposition \eqref{regularity:tmp2} can be rewritten\footnote{We slightly abuse notation and write $\nabla\psi$ for the piecewise gradient of $\psi\in\PSHM{\ell+2}(\Omega)$.} as 
	\begin{align}\label{regularity:tmp3}
		\solv = \solz+\nabla\psi.
	\end{align}
	
	Note that we are not interested in $\psi$ itself, but only in its piecewise gradient $\nabla\psi$. This means that we can alter $\psi$ by adding piecewise constants, hence, without loss of generality we may assume that $\psi$ satisfies $\SCP{\psi}{1}{\SL(\Omega)}=0$ and $\SCP{\jump{\psi}}{1}{\SL(\interf)}=0$.
	
	We notice that \eqref{regularity:tmp3} implies
	$-\diverg\nu\nabla \psi = \diverg(\nu\solz)$ on every subdomain $\Gp_i$, and due to $\solv\in\Hcurl$ and $\diverg \nu\solv=0$ we have $\tjump{\nabla \psi} = -\tjump{\solz}$ and $\njump{\nu\nabla\psi} = -\njump{\nu\solz}$ on $\interf$. Furthermore, the boundary condition $\nu\solv\cdot\soln = 0$ leads to $\nu\nabla\psi\cdot\soln = -\nu\solz\cdot\soln$ on $\Gamma$.

	Hence, according to Lemma \ref{regularity:sgtransmission} there holds $\psi\in\PSHM{\ell+3}(\Gp)$ and 
	\begin{align}\label{reg:tmp4}
	\begin{split}
	\norm{\nabla\psi}{\PVSHM{\ell+2}(\Gp)}&\leq C\left(\norm{\solz}{\PVSHM{\ell+2}(\Gp)}+\norm{\nu\solz\cdot\soln}{\SHM{\ell+3/2}(\Gamma)}+\norm{\tjump{\solz}}{\VSHM{\ell+3/2}_T(\interf)}+\norm{\njump{\nu\solz}}{\SHM{\ell+3/2}(\interf)}\right) \\
	&\leq C\norm{\solz}{\PVSHM{\ell+2}(\Gp)},
	\end{split}
	\end{align}
	where the last estimate comes from the trace inequality applied to the terms living on $\Gamma$ and $\interf$.
	
Together with \eqref{reg:zinequality} and $\eqref{regularity:tmp3}$ this shows
\begin{align*}
	\norm{\solv}{\PVSHM{\ell+2}(\Gp)}\leq C\norm{\curl \solv}{\PVSHM{\ell+1}(\Gp)},
\end{align*}
which finishes the proof of \eqref{sregtmp}.

If we have that $\nu\in\VCM{\ell+1}(\overline{\Omega})$ and $\solu\in\VHcurl{\ell}\cap\VHdiv{\ell}$, the proof follows the same lines but instead of exploiting Lemma~\ref{regularity:sgtransmission}, one can employ the regularity shift properties of the standard Poisson problem with Neumann boundary conditions.
\end{fatproofmod}

Having proved Theorem~\ref{Mainresult1}, we can use it to give a proof for Theorem~\ref{Mainresult2}.

\medskip 

%
%
%

\begin{fatproofmod}{Theorem~\ref{Mainresult2}}
	We notice that due to $\Omega$ being a $\CM{\ell+2}$-domain, there exist constants $c_1,c_2>0$ depending only on $\Omega$ and $\ell$ such that 
	\begin{align*}
	c_1\norm{\solu_T}{\VSHM{\ell+1/2}_T(\Gamma)}\leq \norm{\solu_t}{\VSHM{\ell+1/2}_T(\Gamma)}\leq c_2\norm{\solu_T}{\VSHM{\ell+1/2}_T(\Gamma)}.
\end{align*}
Hence, $\solu_T=\solg_T$ on $\Gamma$ for a tangent field $\solg_T\in\VSHM{\ell+1/2}_T(\Gamma)$ if and only if $\solu_t = \solh_T$ on $\Gamma$ for a tangent field $\solh\in\VSHM{\ell+1/2}_T(\Gamma)$, and without loss of generality we may assume that $\solu$ satisfies $\solu_T=\solg_T$ on $\Gamma$ for a tangent field $\solg_T\in\VSHM{\ell+1/2}_T(\Gamma)$.

	According to Lemma \ref{reg:helmholtzdiv}, there exists a decomposition $\nu\solu = \nu\solv+\curl\solz$ with
	\begin{align*}
			\norm{\solv}{\PVSHM{1}\revision{(\Gp)}}\leq C\norm{\diverg\nu\solu}{\SL(\Omega)},\quad \norm{\solz}{\VSHM{1}(\Omega)}\leq \norm{\solu}{\VSL(\Omega)},
	\end{align*}
	and $\SCP{\solv}{\curl\solw}{\VSL(\Omega)} = 0$ for all $\solw\in\Hcurl$. Thus, the vector field $\solj := \nu^{-1}\curl\solz$ satisfies the equations $\curl\solj = \curl\solu\in\PVSHM{\ell}(\Gp)$, $\diverg\nu\solj = 0$ in $\Omega$ and $\solj_T = \solg$ on $\Gamma$.
	The remainder of the proof is split in two steps.
	
	{\bf Step 1:} The first step is to prove the inequality
	\begin{align}\label{tmptimes0}
		\norm{\solj}{\PVSHM{\ell+1}(\Gp)}\leq C\norm{\curl\solu}{\PVSHM{\ell}(\Gp)}.
	\end{align}
	We employ Lemma \ref{reg:helmholtze} to obtain a splitting $\solj = \solr+\nabla\phi$ with 
	\begin{align*}
			\norm{\solr}{\PVSHM{1}(\Gp)}\leq C\norm{\curl\solu}{\VSL(\Omega)},\quad \norm{\revision{\phi}}{\SHM{1}(\Omega)}\leq C\norm{\solu}{\VSL(\Omega)},
	\end{align*}
	and $\SCP{\nu\solr}{\nabla\xi}{\VSL(\Omega)} = 0$ for all $\xi\in\SHM{1}(\Omega)$. We notice that $\solr$ satisfies $\diverg\nu\solr=0$ in $\Omega$, $\curl\solr = \curl \solj = \curl\solu$ in $\Omega$ and $\nu\solr\cdot\soln = 0$ on $\Gamma$. According to Theorem \ref{Mainresult1}, this implies 
	\begin{align*}
		\norm{\solr}{\PVSHM{\ell+1}(\Gp)}\leq C\norm{\curl\solu}{\PVSHM{\ell}(\Gp)}.
	\end{align*}
	Moreover, there holds $\diverg\nu\nabla\phi = \diverg\nu\solj = 0$ in $\Omega$ and $\nabla_{\Gamma}\phi = \solg-\solr_T\in\VSHM{\ell+1/2}_T(\Gamma)$. Thus, $\diverg_{\Gamma}\nabla_{\Gamma}\phi\in\SHM{\ell-1/2}(\Gamma)$ and consequently, due to elliptic regularity, $\phi\vert_{\Gamma}\in\SHM{\ell+3/2}(\Gamma)$. Therefore, the elliptic regularity of the Poisson problem shows $\norm{\nabla\phi}{\PSHM{\ell+1}(\Gp)}\leq C\left(\norm{\solr}{\PVSHM{\ell+1}(\Gp)}+\norm{\solg}{\VSHM{\ell+1/2}_T(\Gamma)}\right)\leq C\left(\norm{\curl\solu}{\PVSHM{\ell}(\Gp)}+\norm{\solg}{\VSHM{\ell+1/2}_T(\Gamma)}\right)$, which finishes the proof of \eqref{tmptimes0}.
	
	{\bf Step 2:} The second step is to prove 
	\begin{align}\label{tmptimes1}
		\norm{\solv}{\PVSHM{\ell+1}(\Gp)}\leq C\norm{\diverg\nu\solu}{\PVSHM{\ell}(\Gp)}.
	\end{align}
	The property $\SCP{\solv}{\curl\solw}{\VSL(\Omega)} = 0$ for all $\solw\in\Hcurl$ implies that $\solv = \nabla\xi$ for a function $\xi\in\SHM{1}_0(\Omega)$, see for example \cite[Lemma~3.27,~(3.60)]{BookMonk}. Moreover, $\xi$ satisfies $\diverg\nu\nabla\xi = \diverg\nu\solu\in\PVSHM{\ell}(\Gp)$, which, according to Lemma \ref{regularity:sgtransmissiondir} yields $\norm{\xi}{\PSHM{\ell+2}(\Gp)}\leq C\norm{\diverg\nu\solu}{\PVSHM{\ell}(\Gp)}$. This proves \eqref{tmptimes1} and finishes the proof in the case of (possibly) discontinuous $\nu$.

	If we have that $\nu\in\VCM{\ell+1}(\overline{\Omega})$ and $\solu\in\VHcurl{\ell}\cap\VHdiv{\ell}$, the proof follows the same lines but all transmission problems become standard Poisson problems.
\end{fatproofmod}

%% file: main.bbl
\providecommand{\bysame}{\leavevmode\hbox to3em{\hrulefill}\thinspace}
\providecommand{\MR}{\relax\ifhmode\unskip\space\fi MR }
\providecommand{\MRhref}[2]{%
  \href{http://www.ams.org/mathscinet-getitem?mr=#1}{#2}
}
\providecommand{\href}[2]{#2}
\begin{thebibliography}{10}

\bibitem{DecompositionAmrouche}
C.~Amrouche, C.~Bernardi, M.~Dauge, and V.~Girault, \emph{Vector potentials in
  three-dimensional nonsmooth domains}, Math. Meth. Appl. Sci. \textbf{21}
  (1998), 823--864.

\bibitem{TracesHcurl}
A.~Buffa, M.~Costabel, and D.~Sheen, \emph{On traces for {{\bf H}}{\bf
  (curl},{$ \Omega$}{\bf )} in {L}ipschitz domains}, J. Math. Anal. Appl.
  \textbf{276} (2002), no.~2, 845--867.

\bibitem{TransmissionCaffarelli}
L.A. Caffarelli, M.~Soria-Carro, and P.~R. Stinga, \emph{Regularity for
  {$C^{1,\alpha}$} interface transmission problems}, Arch. Ration. Mech. Anal.
  \textbf{240} (2021), no.~1, 265--294.

\bibitem{BookCessenat}
M.~Cessenat, \emph{Mathematical methods in electromagnetism}, World Scientific,
  Singapore, 1996.

\bibitem{MaxwellSpence}
T.~Chaumont-Frelet, J.~Galkowski, and E.A. Spence, \emph{Sharp error bounds for
  edge-element discretisations of the high-frequency {M}axwell equations},
  (2024), Preprint, https://arxiv.org/abs/2408.04507.

\bibitem{MaxwellChaumont}
T.~Chaumont-Frelet and P.~Vega, \emph{Frequency-explicit approximability
  estimates for time-harmonic {M}axwell's equations}, Calcolo \textbf{59}
  (2022), no.~2, Paper No. 22, 15 pp.

\bibitem{MaxwellChen}
Z.~Chen, \emph{On the regularity of time-harmonic {M}axwell equations with
  impedance boundary conditions}, Commun. Appl. Math. Comput. (2024).

\bibitem{MaxwellCostabel}
M.~Costabel, \emph{A remark on the regularity of solutions of {M}axwell's
  equations on {L}ipschitz domains}, Math. Methods Appl. Sci. \textbf{12}
  (1990), no.~4, 365--368.

\bibitem{GrandLivre}
M.~Costabel, M.~Dauge, and S.~Nicaise, \emph{{Corner Singularities and Analytic
  Regularity for Linear Elliptic Systems. Part I: Smooth domains.}}, Preprint,
  https://hal.science/hal-00453934v2, 2010.

\bibitem{CurlInverse}
M.~Costabel and A.~Mc{I}ntosh, \emph{On {B}ogovski\u{\i} and regularized
  {P}oincar\'{e} integral operators for de {R}ham complexes on {L}ipschitz
  domains}, Math. Z. \textbf{265} (2010), no.~2, 297--320.

\bibitem{MaxwellHiptmair}
R.~Hiptmair, A.~Moiola, and I.~Perugia, \emph{Stability results for the
  time-harmonic {M}axwell equations with impedance boundary conditions}, Math.
  Models Methods Appl. Sci. \textbf{21} (2011), no.~11, 2263--2287.

\bibitem{TransmissionNirenberg}
Y.~Li and L.~Nirenberg, \emph{Estimates for elliptic systems from composite
  material}, Comm. Pure Appl. Math. \textbf{56} (2003), no.~7, 892--925,
  Dedicated to the memory of J\"urgen K. Moser.

\bibitem{TransmissionVogelius}
Y.~Li and M.~Vogelius, \emph{Gradient estimates for solutions to divergence
  form elliptic equations with discontinuous coefficients}, Arch. Ration. Mech.
  Anal. \textbf{153} (2000), no.~2, 91--151.

\bibitem{MaxwellLu}
P.~Lu, Y.~Wang, and X.~Xu, \emph{Regularity results for the time-harmonic
  {M}axwell equations with impedance boundary condition},  (2018), Preprint,
  https://arxiv.org/abs/1804.07856.

\bibitem{BookMcLean}
W.~McLean, \emph{Strongly elliptic systems and boundary integral equations},
  Cambridge University Press, Cambridge, 2000.

\bibitem{MaxwellTransparentMelenk}
J.M. Melenk and S.A. Sauter, \emph{Wavenumber-explicit {$hp$}-{FEM} analysis
  for {M}axwell's equations with transparent boundary conditions}, Found.
  Comput. Math. \textbf{21} (2021), no.~1, 125--241.

\bibitem{MaxwellImpedanceMelenk}
\bysame, \emph{Wavenumber-{E}xplicit {\it hp}-{FEM} {A}nalysis for {M}axwell's
  {E}quations with {I}mpedance {B}oundary {C}onditions}, Found. Comput. Math.
  \textbf{24} (2024), no.~6, 1871--1939.

\bibitem{MaxwellMyself2}
J.M. Melenk and D.~Wörgötter, \emph{Wavenumber-explicit analytic regularity
  of {M}axwell's equations with impedance boundary conditions in piecewise
  smooth media}, In preparation, 2025.

\bibitem{MaxwellMyself3}
\bysame, \emph{Wavenumber-explicit {$hp$}-{FEM} analysis of {M}axwell's
  equations with impedance boundary conditions in piecewise smooth media}, In
  preparation, 2025.

\bibitem{BookMonk}
P.~Monk, \emph{Finite element methods for {M}axwell's equations}, Oxford
  {U}niversity {P}ress, {N}ew {Y}ork, 2003.

\bibitem{BookNedelec}
J.C. N\'{e}d\'{e}lec, \emph{Acoustic and {E}lectromegnetic {E}quations},
  Springer, {N}ew {Y}ork, 2001.

\bibitem{MaxwellTomezyk}
S.~Nicaise and J.~Tomezyk, \emph{Convergence analysis of a {$hp$}-finite
  element approximation of the time-harmonic {M}axwell equations with impedance
  boundary conditions in domains with an analytic boundary}, Numer. Methods
  Partial Differential Equations \textbf{36} (2020), no.~6, 1868--1903.

\bibitem{RegularityWeber}
C.~Weber, \emph{Regularity {T}heorems for {M}axwell's {E}quations}, Math. Meth.
  in the Appl. Sci. \textbf{3} (1981), 523--536.

\end{thebibliography}
